\pdfoutput=1
\RequirePackage{fix-cm}
\documentclass[envcountsect]{svjour3}
\smartqed
\usepackage[utf8]{inputenc}
\usepackage{amsmath,amsfonts}
\usepackage{tikz}
\usepackage{graphicx}
\usepackage[capitalize]{cleveref}
\usepackage{mathptmx}

\newcommand{\diag}{\operatorname{diag}\nolimits}
\newcommand{\arcosh}{\operatorname{arcosh}\nolimits}
\newcommand{\erfcx}{\operatorname{erfcx}\nolimits}
\title{Numerical Evaluation of Mittag-Leffler Functions}
\titlerunning{Mittag-Leffler Functions}
\author{William McLean}
\institute{William McLean\at
School of Mathematics and Statistics, University of New South Wales,
Sydney 2052, Australia\\
\email{w.mclean@unsw.edu.au}
}
\date{\today}
\begin{document}
\maketitle
\begin{abstract}
The Mittag-Leffler function is computed via a quadrature approximation of a 
contour integral representation.  We compare results for parabolic and 
hyperbolic contours, and give special attention to evaluation on the real line. 
The main point of difference with respect to similar approaches from the 
literature is the way that poles in the integrand are handled.  Rational 
approximation of the Mittag-Leffler function on the negative real axis is also 
discussed.
\end{abstract}
\section{Introduction}
The Mittag-Leffler function~\cite{GorenfloEtAl2014,MittagLeffler1903} is 
defined by the power series
\begin{equation}\label{eq: E defn}
E_\alpha(z)=\sum_{n=0}^\infty\frac{z^n}{\Gamma(1+n\alpha)},
\end{equation}
which converges for all~$z\in\mathbb{C}$ if~$\alpha>0$. For integer 
values of~$\alpha$ we can express $E_\alpha$ in terms of elementary 
functions~\cite[Equation~(1.5)]{Paris2002}, e.g.,
\begin{equation}\label{eq: integer alpha}
E_0(z)=\frac{1}{1-z},\qquad E_1(z)=e^z,\qquad
E_2(z)=\cosh(z^{1/2})=\sum_{n=0}^\infty\frac{z^n}{(2n)!},
\end{equation}
and for $\alpha=1/2$ in terms of the scaled complementary error 
function~\cite[Equation~(2.7)]{HauboldEtAl2011},
\[
E_{1/2}(z)=\erfcx(-z)\quad\text{where}\quad
\erfcx(z)=\frac{2}{\sqrt{\pi}}\int_z^\infty\exp(z^2-t^2)\,dt.
\]
Numerical evaluation of the series~\eqref{eq: E defn} is efficient for 
small~$|z|$, but other approaches are required in general.

The Mittag--Leffler function has been mostly neglected by those responsible for 
mathematical software libraries, apparently because until the 1990s there were 
few applications that called for numerical values of~$E_\alpha(z)$. 
Mathematica seems to provide the only officially supported implementation with 
its \texttt{MittagLefflerE} function.  Growing interest in the applications of 
fractional calculus, such as in fractional partial differential equation 
models, led in the 2000s to a growing demand for methods to 
evaluate~$E_\alpha(z)$, particularly for the case~$0<\alpha<1$ with $z$ on the 
negative real axis.  Such applications involve also the two-parameter 
Mittag-Leffler function,
\begin{equation}\label{eq: gen E defn}
E_{\alpha,\beta}(z)=\sum_{n=0}^\infty\frac{z^n}{\Gamma(\beta+n\alpha)},
\end{equation}
as well as other generalizations beyond the scope of this paper.  The 
identity~\cite[Equation~(5.4)]{HauboldEtAl2011}
\[
\frac{d}{dz}\,E_{\alpha,\beta}(z)=\frac{1}{\alpha z}\bigl[E_{\alpha,\beta-1}(z)
    -(\beta-1)E_{\alpha,\beta}(z)\bigr]
\]
means that if the Mittag-Leffler function can be evaluated then so can 
its derivative.

Gorenflo et al.~\cite{GorenfloEtAl2002} and Seybold and 
Hilfer~\cite{SeyboldHilfer2008} developed methods to 
compute~$E_{\alpha,\beta}(x)$ based on contour integral representations and
asymptotic expansions for large~$|z|$.  These authors provide computable error 
bounds that reduce the problem to evaluating integrals of complex-valued 
functions defined on finite real intervals, which then require some kind of 
numerical quadrature. Igor Podlubny~\cite{Podlubny_mlf} wrote a widely used, 
third-party Matlab function, \texttt{mlf}, based on these methods.

Weideman and Trefethen~\cite{WeidemanTrefethen2007} developed a quadrature 
method for numerical inversion of the Laplace transform via the Bromich 
integral formula, thereby providing another way to evaluate~$E_\alpha$.  
The method was developed further by  
Garrappa~\cite{Garrappa2015,GarrappaPopolizio2013} who contributed another 
third-party Matlab code, \texttt{ml}.  Gill and Straka~\cite{GillStraka2017}
used a similar approach in their R package, \texttt{MittagLeffleR}, which
handles Mittag-Leffler probability distributions.  

The approach used here is similar to that of Garrappa~\cite{Garrappa2015}, but 
whereas he used the contour integral representation
\begin{equation}\label{eq: Garrappa integral}
e_{\alpha,\beta}(t;\lambda)=t^{\beta-1}E_{\alpha,\beta}(t^\alpha\lambda)
    =\frac{1}{2\pi i}\int_{\mathcal{C}}e^{st}\,
    \frac{s^{\alpha-\beta}\,ds}{s^\alpha-\lambda},
\quad t>0,\quad\lambda\in\mathbb{C},
\end{equation}
we will work with the Wiman integral,
\begin{equation}\label{eq: Wiman integral}
E_{\alpha,\beta}(z)=\frac{1}{2\pi i}\int_{\mathcal{C}}
    \frac{e^w w^{\alpha-\beta}}{w^\alpha-z}\,dw.
\end{equation}
In both cases, one makes a branch cut along the negative real axis (in the 
$s$-plane~and $w$-plane respectively) and chooses a Hankel 
contour~$\mathcal{C}$ that encircles this cut in the counterclockwise direction, 
beginning at infinity in the third quadrant, passing to the right of all 
singularities in the integrand, and finishing at infinity in the second 
quadrant.  The two representations are related via the substitution~$w=st$.

When applying the method of Weideman and Trefethen to 
evaluate~\eqref{eq: Garrappa integral}, the contour~$\mathcal{C}$ is chosen to 
be a parabola or hyperbola that depends on~$t$.  The parameters describing this 
contour, together with the step size in the quadrature formula, are chosen to 
optimize the convergence rate.  We use essentially the same approach 
for~\eqref{eq: Wiman integral} but with a fixed~$\mathcal{C}$ independent 
of~$z$.  The more substantial difference in our approach concerns the treatment 
of any poles that arise when the denominator of the integrand vanishes.  Whereas 
Garrappa avoids the poles by adjusting $\mathcal{C}$, we instead split the 
integrand into a simple singular term and a remainder that is analytic in the 
cut plane.

After summarizing a few standard properties of the Mittag-Leffler function in
\cref{sec: integral rep}, we present the quadrature method in 
\cref{sec: Quadrature}.  If the argument~$z$ is real, then the 
computational cost of the method can be reduced using techniques we describe in 
\cref{sec: real line}.  For comparison, \cref{sec: Pade} considers rational 
approximations of Pad\'e type that were introduced by Zeng and 
Chen~\cite{ZengChen2015}.  Finally, in~\cref{sec: best approx} we present some 
numerical examples of minimax approximation of~$E_\alpha(-x)$ for~$0\le 
x<\infty$ using rational functions of type~$(m,m)$.  

A Julia package~\cite{MLJulia} implements the numerical methods described 
herein.  The scripts used to generate the figures and tables below are included 
in the \texttt{examples} folder of the package repository on Github.  All 
computations were performed in 64-bit floating point arithmetic although some 
routines in the package support the use of Julia's \texttt{BigFloat} multiple 
precision data type.

\section{Integral representation and asymptotics}\label{sec: integral rep}
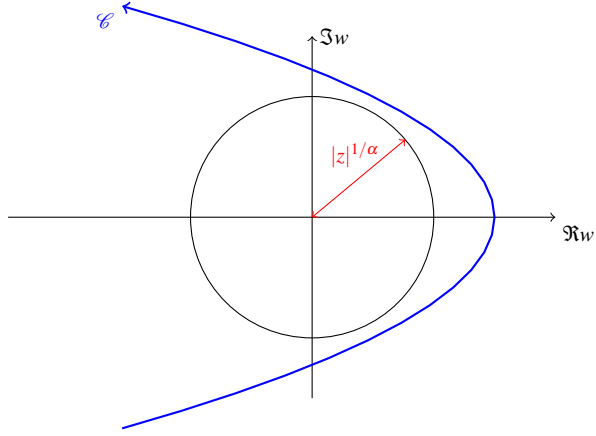
\begin{figure}
\begin{center}
\begin{tikzpicture}[scale=0.8]
\draw[->] (-5,0) -- (4,0);
\node[below right] at (4,0) {$\Re w$};
\draw[->] (0,-3) -- (0,3);
\node[right] at (0,3) {$\Im w$};
\draw[->,thick,domain=-3.5:3.5, blue] plot (3-\x*\x/2, \x);
\draw[-] (0, 0) circle [radius=2];
\node[blue, below left] at ({-25/8}, 3.5) {$\mathcal{C}$};
\draw[<->,red] (0,0) -- ({2*cos(40)}, {2*sin(40)});
\node[red,above] at ({sqrt(2)/2},{sqrt(2)/2}) {$|z|^{1/\alpha}$};
\end{tikzpicture}
\end{center}
\caption{The Hankel contour $\mathcal{C}$.}\label{fig: Hankel}
\end{figure}

Properties of the Mittag-Leffler function 
are described at length in the survey 
article of Haubold et al.~\cite{HauboldEtAl2011} and the mongraph of Gorenflo 
et al.~\cite{GorenfloEtAl2014}.  We will require only a few simple facts
that are easily derived using the integral representation of the reciprocal of 
the Gamma function:
\begin{equation}\label{eq: Gamma Hankel}
\frac{1}{\Gamma(\nu)}=\frac{1}{2\pi i}\int_{\mathcal{C}}\frac{e^w}{w^\nu}\,dw.
\end{equation}
By choosing $\mathcal{C}$ so that it passes to the right of the circle 
$|w|=|z|^{1/\alpha}$, as illustrated in \cref{fig: Hankel}, and using
\eqref{eq: Gamma Hankel} with~$\nu=\beta+n\alpha$ in~\eqref{eq: gen E defn},
the integral representation~\eqref{eq: Wiman integral} follows by summation of 
a geometric series, which converges because $|z/w^\alpha|<1$ for all 
$w\in\mathcal{C}$.  

Suppose now that $0<\alpha<1$ and $z\ne0$ with
$\theta=\arg z\in(-\pi,\pi]$.  If $\alpha\pi<|\theta|<\pi$, 
then $w^\alpha-z\ne0$ for~$w$ in the cut plane $|\arg w|<\pi$ and we find that
$E_{\alpha,\beta}(z)=O(|z|^{-1})$ as~$|z|\to\infty$.  However, if 
$|\theta|<\alpha\pi$ then the equation~$w^\alpha-z=0$ has a single solution
in the cut plane, namely $w=\gamma$ where
\begin{equation}\label{eq: gamma}
\gamma=z^{1/\alpha}=|z|^{1/\alpha}\exp(i\theta/\alpha).
\end{equation}
Shifting the contour so it passes to the left of the pole at~$\gamma$ we 
collect a residue, with the result that
\begin{equation}\label{eq: E gamma}
E_{\alpha,\beta}(z)=\alpha^{-1}\gamma^{1-\beta}\exp(\gamma)
    +\frac{1}{2\pi i}\int_{\mathcal{C}'}
    \frac{e^w w^{\alpha-\beta}}{w^\alpha-z}\,dw,
\end{equation}
and the integral over the shifted contour~$\mathcal{C}'$ is again $O(|z|^{-1})$ 
as~$|z|\to\infty$.

The identity
\begin{equation}\label{eq: beta-alpha identity}
E_{\alpha,\beta-\alpha}(z)=\frac{1}{\Gamma(\beta-\alpha)}+zE_{\alpha,\beta}(z)
\end{equation}
follows easily from the power series~\eqref{eq: gen E defn}, leading to
the asymptotic formulae
\begin{align*}
E_{\alpha,\beta}(z)&=\frac{-z}{\Gamma(\beta-\alpha)}
    +z^{-1}E_{\alpha,\beta-\alpha}(z)\\
    &=-\frac{z^{-1}}{\Gamma(\beta-\alpha)}+\begin{cases}
    O(|z|^{-2})&\text{if $\alpha\pi<|\theta|\le\pi$,}\\
    \alpha^{-1}z^{(1-\beta)/\alpha}\exp(z^{1/\alpha})+O(|z|^{-2})
    &\text{if $|\theta|<\alpha\pi$,}
\end{cases}
\end{align*}
Iterating this result, we find that~\cite[Equation (1.2)]{Paris2020}
\begin{equation}\label{eq: asymptotics}
E_{\alpha,\beta}(z)= 
    -\sum_{n=1}^{m-1}\frac{z^{-n}}{\Gamma(\beta-n\alpha)}+\begin{cases}
    O(|z|^{-m})&\text{if $\alpha\pi<|\theta|\le\pi$,}\\
    \alpha^{-1}z^{(1-\beta)/\alpha}\exp(z^{1/\alpha})+O(|z|^{-m})
    &\text{if $|\theta|<\alpha\pi$,}\\
\end{cases}
\end{equation}
The behaviour around the Stokes lines $\theta=\pm\alpha\pi$ is quite subtle and 
has been studied in detail by Paris~\cite{Paris2020,Paris2002} and by Wong and 
Zhao~\cite{WongZhao2002}.

\begin{figure}
\begin{center}
\includegraphics[scale=0.6]{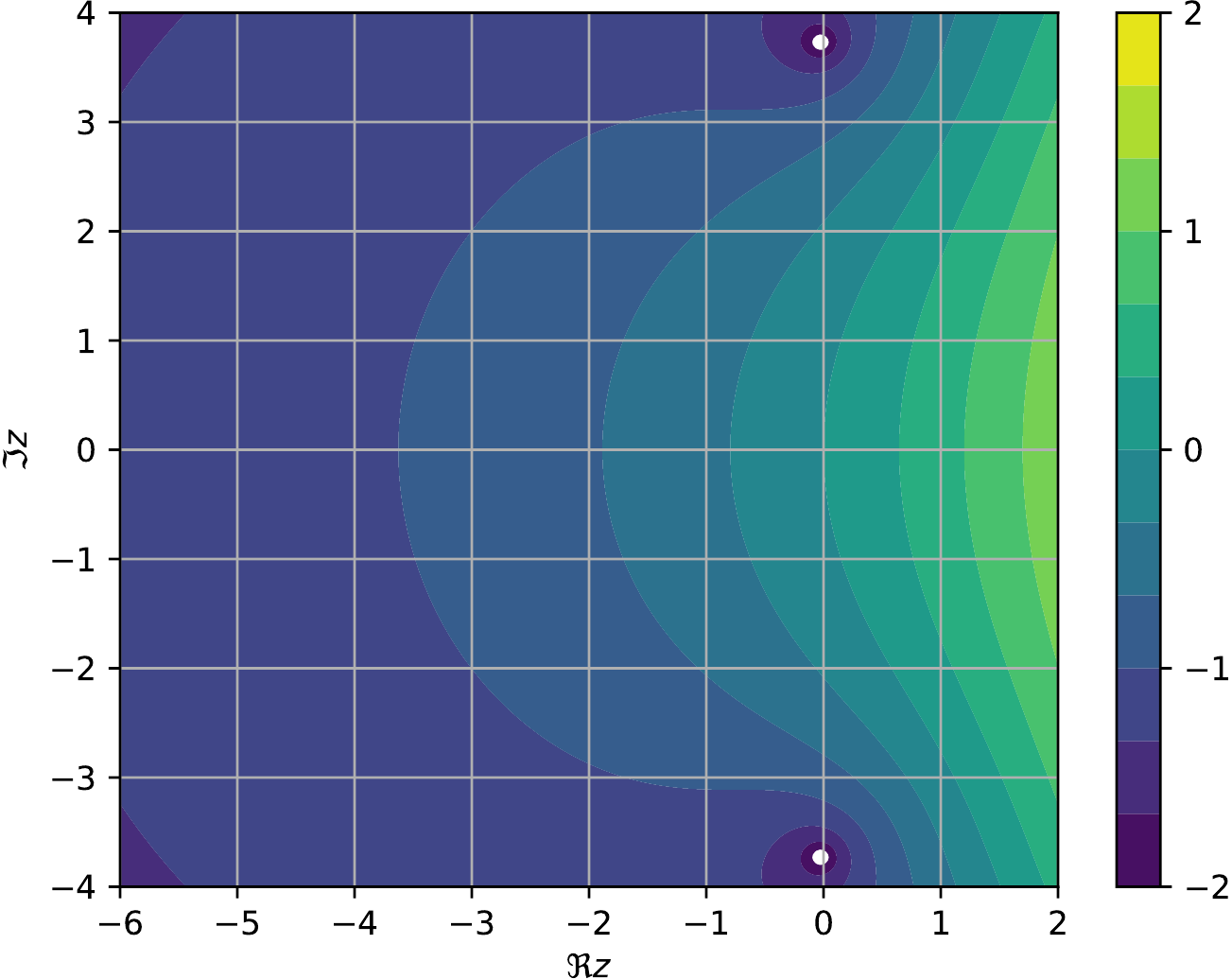} 
\end{center}
\caption{Contour plot of $\log_{10}|E_{3/4}(z)|$.}\label{fig: E contours}
\end{figure}

\begin{example}
Since $|\exp(z^{1/\alpha})|=\exp(|z|^{1/\alpha}\cos\theta/\alpha)$
and $\cos\theta/\alpha>0$ for~$|\theta|<\alpha\pi/2$, it follows that in this 
sector $E_{\alpha,\beta}(z)$ exhibits super exponential growth.  This 
behaviour can be seen, for $\alpha=3/4$, in \cref{fig: E contours}, which shows 
a contour plot of $\log_{10}|E_{3/4}(z)|$ for~$-6\le\Re z\le2$~and 
$-4\le\Im z\le4$.  Notice the zeros near $\pm3{\cdot}7i$.
\end{example}

\begin{remark}\label{remark: asymp}
When $|z|$ is sufficiently large, dropping the $O(|z|^{-m})$ remainder term 
from the asymptotic formula~\eqref{eq: asymptotics} yields a 
practical numerical approximation
$E_{\alpha,\beta}(z)\approx E_{\alpha,\beta,m}^{\mathrm{asymp}}(z)$.
Recalling the identity
\begin{equation}\label{eq: Gamma reflection}
\Gamma(z)\Gamma(1-z)=\frac{\pi}{\sin\pi z},
\end{equation}
we write the coefficient in the $n$th term of the asymptotic sum as 
$1/\Gamma(\beta-n\alpha)=\sigma_n\tau_n$ where
\[
\sigma_n=1\quad\text{and}\quad\tau_n=\frac{1}{\Gamma(\beta-n\alpha)}
\quad\text{if $n\alpha<\beta$,}
\]
with
\[
\sigma_n=-\sin\pi(n\alpha-\beta)\quad\text{and}\quad
\tau_n=\frac{1}{\pi}\,\Gamma(1+n\alpha-\beta)
\quad\text{if $n\alpha\ge\beta$.}
\]
The accuracy of~$E_{\alpha,\beta,m}^{\mathrm{asymp}}(z)$ will improve with 
increasing~$m$ until $m\approx\alpha^{-1}|z|^{1/\alpha}$ 
\cite[Section~2]{Paris2002}, after which the asymptotic sum begins to diverge.
Thus, for a given error tolerance~$\mathrm{tol}$, a simple heuristic is to keep 
adding terms until $\tau_n|z|^{-n}<\mathrm{tol}$ or 
$n>\alpha^{-1}|z|^{1/\alpha}$.  If the former occurs first, then we can expect 
to achieve the desired accuracy; in either case, the size of 
$\tau_{m-1}|z|^{-(m-1)}$ should give an indication of the error.  For a more 
rigorous approach to estimating the remainder in the asymptotic expansion, see 
Seybold and Hilfer~\cite[Theorem~4.2]{SeyboldHilfer2008}.
\end{remark}

\begin{table}
\caption{Results approximating $E_{\alpha,\beta}(-x)$ by applying the algorithm 
of \cref{remark: asymp}, based on the asymptotic 
expansion~\eqref{eq: asymptotics}, when $\alpha=0.7$, $\beta=1.0$ and 
$\mathrm{tol}=10^{-12}$.}
\label{table: asymp}
\begin{center}
\renewcommand{\arraystretch}{1.2}
\begin{tabular}{r|r|r|r|c|c}
$x$&$m$&$\alpha^{-1}x^{1/\alpha}$&\multicolumn{1}{c|}{Error}&
$\tau_{m-1}|z|^{-(m-1)}$&$\tau_m|z|^{-m}$\\ 
\hline
     5&    15&     14.2&   5.84e-06&   1.21e-05&   1.18e-05\\
    15&    16&     68.4&   8.11e-14&   8.24e-13&   2.82e-13\\
    25&    12&    141.9&  -6.54e-14&   3.70e-13&   6.09e-14\\
    35&    10&    229.5&  -1.25e-14&   8.15e-13&   8.31e-14\\
    45&    10&    328.6&  -4.09e-15&   8.49e-14&   6.73e-15\\
    55&     9&    437.7&   8.24e-15&   2.34e-13&   1.39e-14
\end{tabular}
\end{center}
\end{table}

\begin{example}\label{example: asymp}
Let $z=-x$ for~$x>0$ and consider the error 
$E_{\alpha,\beta}(-x)-E_{\alpha,\beta,m}^{\mathrm{asymp}}(-x)$.
\cref{table: asymp} shows some results applying the method of 
\cref{remark: asymp} with $\alpha=0.7$, $\beta=1$~and $\mathrm{tol}=10^{-12}$.
We see that for~$x=5$ the asymptotic series is at best accurate only to about 
$5$~decimal places, but for all other cases the error is smaller than 
$\mathrm{tol}$ for values of~$m$ that are much smaller than 
$\alpha^{-1}x^{1/\alpha}$.  In all cases, the error is smaller than 
$\tau_{m-1}x^{-(m-1)}$.
\end{example}

For~$\alpha>1$, the equation~$w^\alpha-z=0$ possesses more 
solutions with~$|\arg w|\le\pi$, complicating the behaviour 
of~$E_{\alpha,\beta}(z)$.  Fortunately, once we are able to evaluate the 
Mittag-Leffler function for~$0<\alpha\le1$, we can handle the case~$\alpha>1$ 
using the identity~\cite[Equation (3.8)]{HauboldEtAl2011} 
\begin{equation}\label{eq: E alpha/m}
E_{\alpha,\beta}(z)=\frac{1}{m}\sum_{k=0}^{m-1}
    E_{\alpha/m,\beta}\bigl(z^{1/m}e^{i2\pi k/m}\bigr)
\quad\text{for $m\in\{1,2,3,\ldots\}$,}
\end{equation}
by choosing $m$ to be the unique integer satisfying~$m-1<\alpha\le m$, so 
that $0<\alpha/m\le1$.  
\section{Quadrature method}\label{sec: Quadrature}
We now consider the numerical evaluation of the integral
\begin{equation}\label{eq: f w z}
E_{\alpha,\beta}(z)=\frac{1}{2\pi i}\int_{\mathcal{C}}e^w f(w;z)\,dw
\quad\text{where}\quad f(w;z)=\frac{w^{\alpha-\beta}}{w^\alpha-z},
\end{equation}
assuming $0<\alpha<1$ and, to begin with, $\alpha\pi<|\theta|\le\pi$ so 
$f(w;z)$ is analytic for all~$w$ in the cut plane.

\subsection{Parabolic contour}
We adapt the approach of Weideman and Trefethen~\cite{WeidemanTrefethen2007},
who proposed a Hankel contour $\mathcal{C}=\{\,w(u):-\infty<u<\infty\,\}$ where
\begin{equation}\label{eq: parabola}
w(u)=\mu(1+iu)^2\quad\text{for $-\infty<u<\infty$, with $\mu>0$.}
\end{equation}
Since 
\begin{equation}\label{eq: Re Im w par}
\Re w(u)=\mu(1-u^2)\quad\text{and}\quad\Im w(u)=2\mu u,
\end{equation}
we see that $\mathcal{C}$ is a parabola cutting the real axis at~$\mu$, 
with~$\Re w(u)\to-\infty$ as $|u|\to\infty$.  Write
\begin{equation}\label{eq: f g}
\frac{1}{2\pi i}\int_{\mathcal{C}}e^w f(w;z)\,dw
    =\int_{-\infty}^\infty g(u;z)\,du\quad\text{where}\quad
g(u;z)=\frac{e^{w(u)}}{2\pi i}\,f\bigl(w(u);z\bigr)w'(u)
\end{equation}
and, for a suitable step size~$h>0$, let us seek to approximate this integral by 
an infinite sum
\[
Q_h(f;z)=h\sum_{n=-\infty}^\infty g(nh;z).
\]

The error analysis for~$Q_h(f;z)$ begins by extending the parametric 
representation~\eqref{eq: parabola} to a conformal 
mapping~$\zeta=u+iv\mapsto w(u+iv)=\mu(1-v+iu)^2$.  Since
\[
\Re w(u+iv)=\mu\bigl((1-v)^2-u^2\bigr)\quad\text{and}\quad
\Im w(u+iv)=2\mu(1-v)u,
\]
we see that $u\mapsto w(u+iv)$ is again a parabola with $\Re w(u+iv)\to-\infty$ 
as~$|u|\to\infty$.  It is necessary to assume $v<1$ so that $\Im w(u+iv)$ is an 
increasing function of~$u$.  In fact, as $v$ increases to~$1$ the parabola 
collapses onto the cut along the negative real axis.  Thus, if $0<r<1$~and 
$s>0$, then $g(\zeta;z)$ is analytic on the strip~$-s\le\Im\zeta\le r$ and we 
have the error bound~\cite[Theorem~5.2]{McNameeEtAl1971}
\begin{equation}\label{eq: Qh error}
\biggl|Q_h(f,z)-\int_{-\infty}^\infty g(u;z)\,du\biggr|
    \le\frac{M(r;z)}{\exp(2\pi r/h)-1}
    +\frac{M(-s;z)}{\exp(2\pi s/h)-1},
\end{equation}
where $M(v;z)=\int_{-\infty}^\infty|g(u+iv;z)|\,du$.

In practice, we must truncate the infinite sum and compute, for some positive
integer~$N$,
\begin{equation}\label{eq: QhN}
Q_{h,N}(f;z)=h\sum_{n=-N}^N g(nh;z),
\end{equation}
which, in view of~\eqref{eq: Re Im w par}, leads to an additional error of 
order~$\exp\bigl(\mu(1-(Nh)^2)\bigr)$ from the sum over~$|n|\ge N+1$.  Putting 
$r=r_\delta=1-\delta$, we have
\[
M(r_\delta;z)\le C_\delta\exp\bigl(\mu\delta^2\bigr)\quad\text{and}\quad
M(s;z)\le C_s\exp\bigl(\mu(1+s)^2\bigr),
\]
and hence the overall quadrature error~$\varepsilon_N=\varepsilon_N(\mu,r,s,h)$ 
is of order
\[
\exp(\mu\delta^2-2\pi r_\delta/h)+\exp\bigl(\mu(1+s)^2-2\pi s/h\bigr)
+\exp\bigl(\mu(1-(Nh)^2)\bigr).
\]
Only the second term depends on~$s$, so we minimise
\[
\mu(1+s)^2-\frac{2\pi s}{h}=\mu\biggl[\biggl(s+1-\frac{2\pi}{\mu h}\biggr)^2
    +\frac{2\pi}{\mu h}-\biggl(\frac{\pi}{\mu h}\biggr)^2\biggr]
\]
by choosing $s+1=2\pi/(\mu h)$.  To balance the three error terms we then 
require
\[
\mu\delta^2-\frac{2\pi r_\delta}{h}=\frac{2\pi}{h}-\frac{\pi^2}{\mu h^2}
=\mu\bigl(1-(Nh)^2\bigr).
\]
Neglecting $\delta$, the optimal parameters are
$h_\star=3N^{-1}$ and $\mu_\star=(\pi/12)\,N$,
giving a quadrature error~$\varepsilon_N(\mu_\star,h_\star)$ of 
order~$\exp(-2\pi/h_\star)=\exp(-2\pi N/3)\approx8{\cdot}12^{-N}$.  

For this optimized choice of the parameters, we will denote the quadrature 
sum by
\begin{equation}\label{eq: QN opt}
Q_{\star,N}(f;z)=Q_{h_\star,N}(f;z),
\end{equation}
Since $w'(u)=2\mu_\star i(1+iu)$ we have
\[
\frac{h_\star w'(u)}{2\pi i}=\frac{\mu_\star h_\star}{\pi}\,(1+iu)
=\frac{1+iu}{4}
\]
so
\begin{equation}\label{eq: QhN Cn}
Q_{\star,N}(f;z)=A\sum_{n=-N}^N C_n f\bigl(w(nh_\star);z\bigr)
\end{equation}
where
\begin{equation}\label{eq: A Cn parabola}
A=\frac{1}{4}\quad\text{and}\quad C_n=e^{w(nh_\star)}\,(1+inh_\star).
\end{equation}
The properties
\begin{equation}\label{eq: conj}
w(-u)=\overline{w(u)}\quad\text{and}\quad C_{-n}=\overline{C}_n,
\end{equation}
mean that it suffices to store the points~$w(nh_\star)$ and coefficients~$C_n$ 
for~$0\le n\le N$.  Note, however, that our notation hides the fact that 
$w(nh_\star)$~and $C_n$ depend not only on~$n$ but also on~$N$.  Thus, 
increasing $N$ means recomputing \emph{all} of the points and coefficients.

\begin{remark}
The error analysis above implicitly assumes that the factor~$e^{w(u)}$ 
dominates the influence of the integrand~$g(u;z)$ in~\eqref{eq: f g}.  On 
the one hand, since $f(w;z)=O(|w|^{\alpha-\beta})$ as~$|w|\to0$, we see that 
the 
factor~$M(r;z)$ in the error bound~\eqref{eq: Qh error} will grow for large 
positive~$\beta$.  On the other hand, since $f(w;z)=O(|w|^{-\beta})$ 
as~$|w|\to\infty$, the factor~$M(-s;z)$ will grow for large negative~$\beta$.
To alleviate the growth of~$M(r;z)$, Garrappa~and 
Popolizio~\cite[Section~3]{GarrappaPopolizio2013} modified the above approach 
by considering a restriction~$r\le\textrm{const}<1$, so that the 
parabola~$u\mapsto w(u+ir)$ does not pass too close to the origin. An 
alternative is to reduce the value of~$\beta$ used in the quadrature 
approximation via the identity~\cite[Section~5]{HauboldEtAl2011}
\[
E_{\alpha,\beta}(z)=z^{-m}E_{\alpha,\beta-m\alpha}(z)-\sum_{n=1}^m
    \frac{z^{-n}}{\Gamma(\beta-n\alpha)}.
\]
Similarly, growth in~$M(-s;z)$ for large negative~$\beta$ may be ameliorated via
the identity
\[
E_{\alpha,\beta}(z)=z^mE_{\alpha,\beta+m\alpha}(z)+\sum_{n=0}^{m-1}
\frac{z^n}{\Gamma(\beta+n\alpha)}.
\]
\end{remark}

\subsection{Hyperbolic contour}
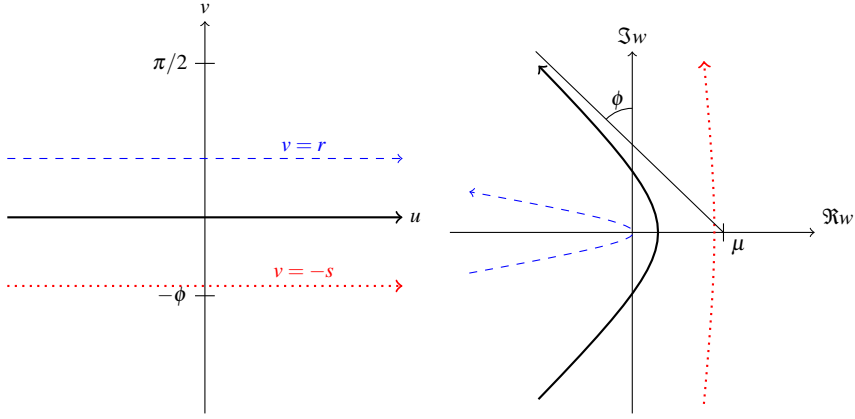
\begin{figure}
\begin{center}
\pgfmathsetmacro{\rval}{0.6}
\pgfmathsetmacro{\sval}{0.7}
\pgfmathsetmacro{\phival}{0.8}
\begin{tikzpicture}[scale=1.3]
\draw[->,thick] (-2,0) -- (2,0);
\node[right] at (2,0) {$u$};
\draw[->] (0,-2) -- (0,2);
\node[above] at (0,2) {$v$};
\draw[->,dashed,blue] (-2,\rval) -- (2,\rval);
\node[above,blue] at (1,\rval) {$v=r$};
\draw[->,red,dotted,thick] (-2,-\sval) -- (2,-\sval);
\node[above,red] at (1,-\sval) {$v=-s$};
\draw[-] (-0.1,-\phival) -- (0.1,-\phival);
\node[left] at (-0.1,-\phival) {$-\phi$};
\draw[-] (-0.1,1.57) -- (0.1,1.57);
\node[left] at (-0.1,1.57) {$\pi/2$};
\end{tikzpicture}
\pgfmathsetmacro{\phideg}{deg(\phival)}
\pgfmathsetmacro{\phipr}{deg(\phival+\rval)}
\pgfmathsetmacro{\phims}{deg(\phival-\sval)}
\begin{tikzpicture}[scale=1.2]
\draw[->] (-2,0) -- (2.0,0);
\node[above right] at (2,0) {$\Re w$};
\draw[->] (0,-2) -- (0,2);
\node[above] at (0,2) {$\Im w$};
\draw[-] (1,-0.1) -- (1,0.1);
\node[below right] at (1,0) {$\mu$};
\draw[->,domain=-1.7:1.7,dashed,blue] plot
    ({1-sin(\phipr)*cosh(\x)},{cos(\phipr)*sinh(\x)});
\draw[->,thick,domain=-1.7:1.7] plot
    ({1-sin(\phideg)*cosh(\x)},{cos(\phideg)*sinh(\x)});
\draw[->,domain=-1.4:1.4,dotted,thick,red] plot
    ({1-sin(\phims)*cosh(\x)},{cos(\phims)*sinh(\x)});
\draw[-,thin] (1,0) -- ({1-2*tan(\phideg)},2);
\draw[-] (0,{0.4+cot(\phideg)}) arc
    [radius=0.4, start angle=90, end angle={90+\phideg}];
\node[above left] at (0,{0.3+cot(\phideg)}) {$\phi$};
\end{tikzpicture}
\end{center}
\caption{Hyperbolas parameterised by $u\mapsto w(u+iv)$.}
\label{fig: hyperbolas}
\end{figure}

The second type of contour considered by Weideman and 
Trefethen~\cite{WeidemanTrefethen2007} is of the form
\begin{equation}\label{eq: hyperbola}
w(u)=\mu\bigl(1+\sin(iu-\phi)\bigr)\quad\text{for $-\infty<u<\infty$,}
\end{equation}
where the parameters $\mu$~and $\phi$ satisfy $\mu>0$ and $0<\phi<\pi/2$;
see also L\'opez-Fern\'andez and Palencia~\cite{LopezFernandezPalencia2004}.
Since $\Re w=\mu(1-\cosh u\,\sin\phi)$~and $\Im w=\mu\sinh u\,\cos\phi$, we have
\[
\biggl(\frac{\Re w-\mu}{\mu\sin\phi}\biggr)^2
    -\biggl(\frac{\Im w}{\mu\cos\phi}\biggr)^2=1,
\]
so $\mathcal{C}$ is the left branch of an hyperbola with asymptotes
$\Im w=\pm(\Re w-\mu)\cot\phi$.  We extend \eqref{eq: hyperbola} to a conformal 
mapping $\zeta=u+iv\mapsto w(u+iv)=\mu\bigl[1+\sin(iu-(\phi+v)\bigr)\bigr]$, 
and see that for a fixed~$v$ with $0<\phi+v<\pi/2$ the curve $u\mapsto w(u+iv)$
is the left branch of an hyperbola with asymptotes
$\Im w=\pm(\Re w-\mu)\cot(\phi+v)$, as illustrated in \cref{fig: hyperbolas}
for $v=r$, $0$~and $-s$.  Noting that
\[
\Re w(u+iv)\le\mu\bigl(1-\sin(\phi+v)\bigr)
\quad\text{for $-\infty<u<\infty$,}
\]
it follows that the error bound~\eqref{eq: Qh error} can be applied for 
$0<r<\pi/2-\phi$~and $0<s<\phi$ with
\[
M(r_\delta;x)\le C\exp\bigl(\mu(1-\cos\delta)\bigr)
\quad\text{and}\quad
M(-s_\delta;x)\le C\exp\bigl(\mu(1-\sin\delta)\bigr),
\]
where $r_\delta=\pi/2-\phi-\delta$~and 
$s_\delta=\phi-\delta$.

\begin{figure}
\begin{center}
\includegraphics[scale=0.4]{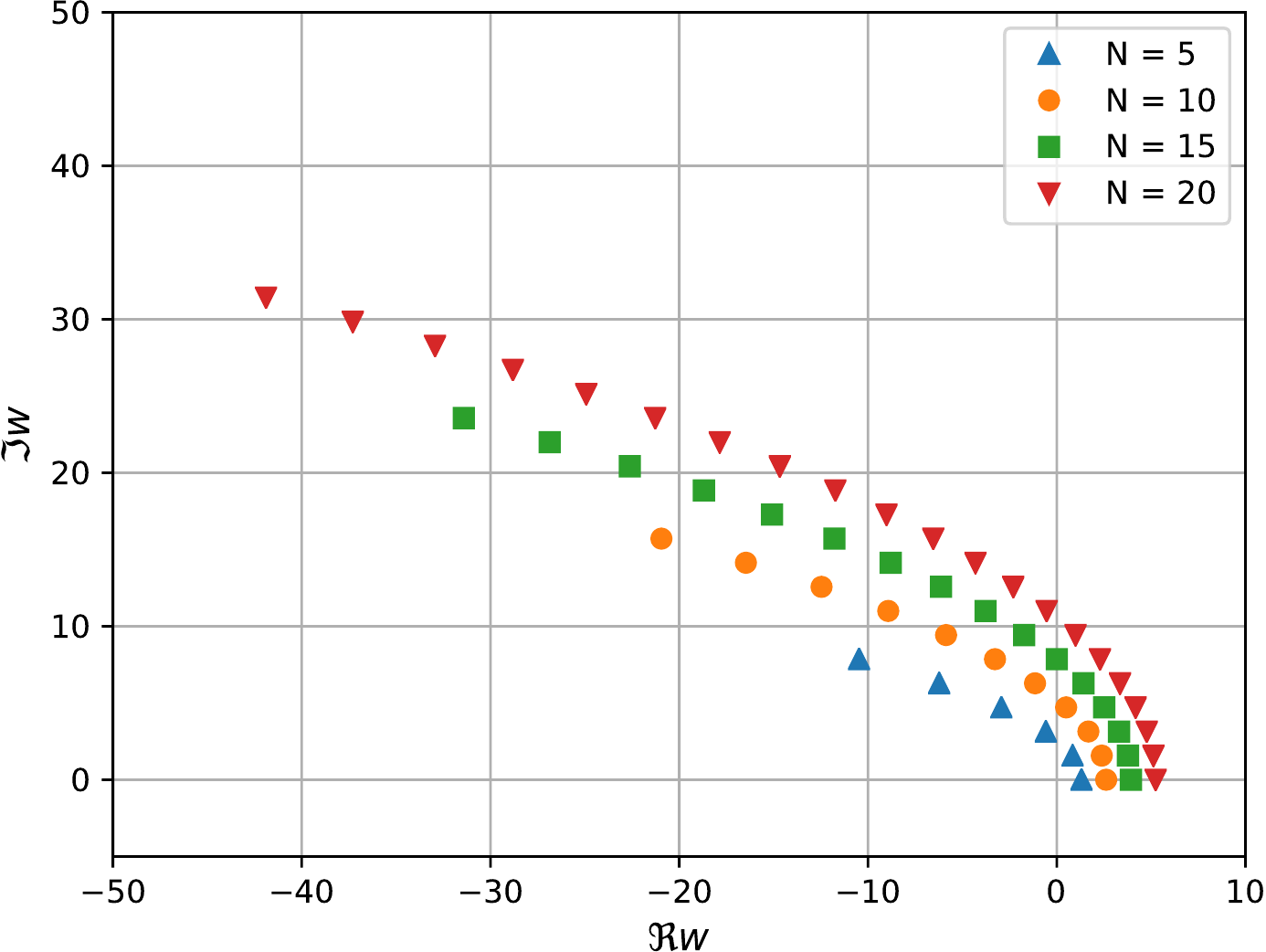}
\hspace{0.2cm}
\includegraphics[scale=0.4]{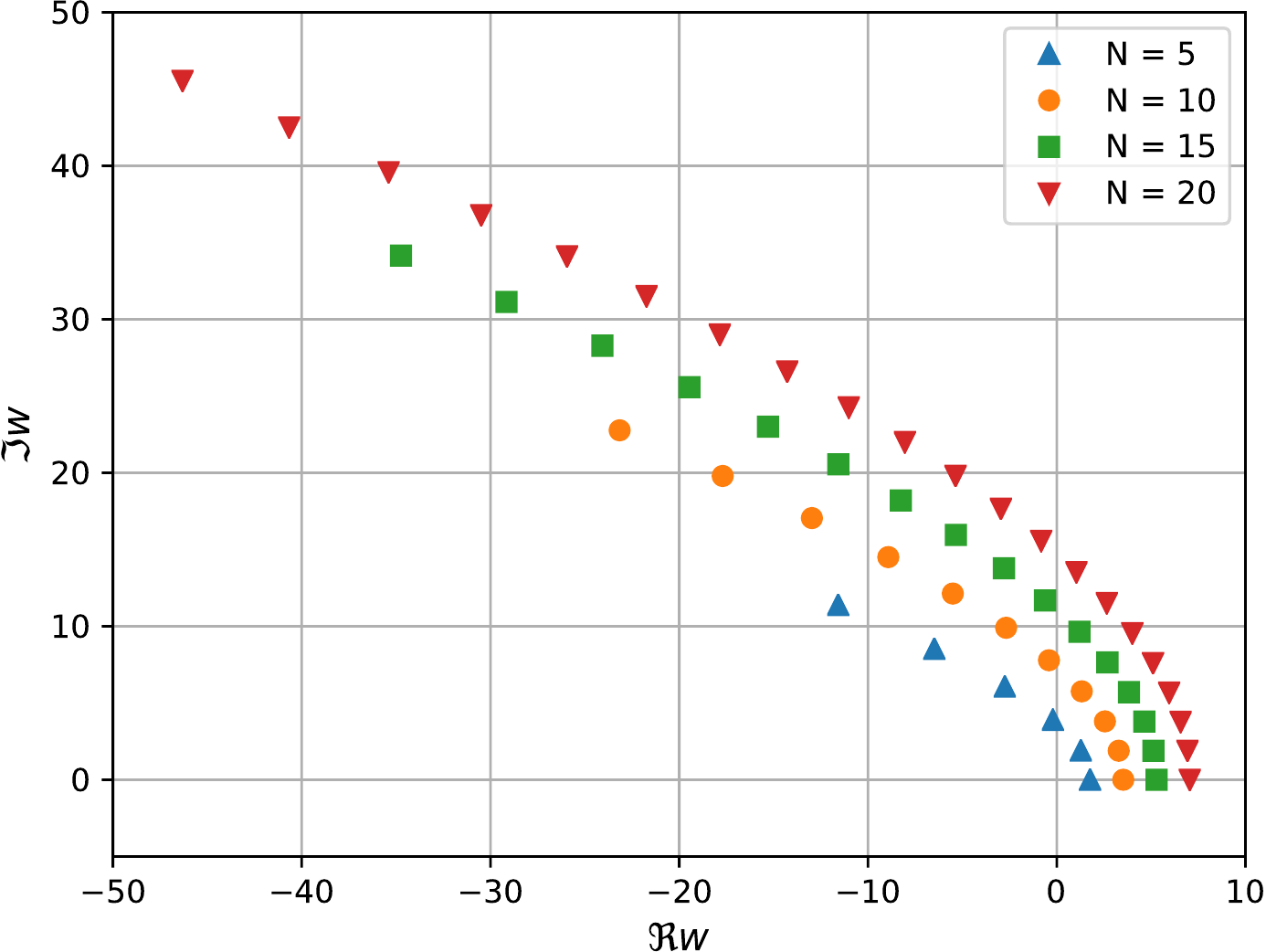}
\end{center}
\caption{The quadrature points $w(nh_\star)$ for $0\le n\le N$ and four choices 
of~$N$ using parabolic (left) and hyperbolic (right) contours.}
\label{fig: quad points}
\end{figure}

\begin{figure}
\begin{center}
\includegraphics[scale=0.4]{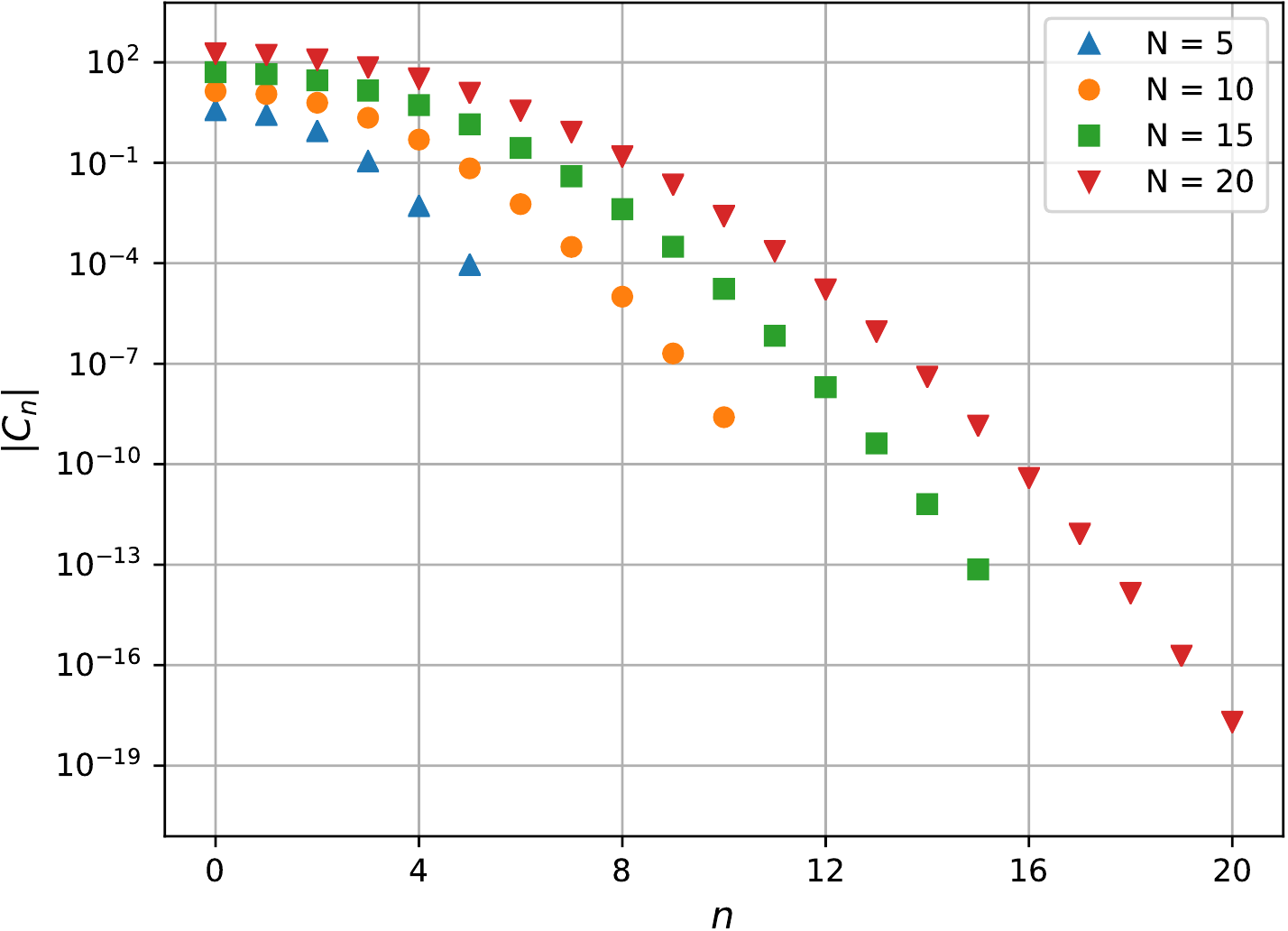}
\hspace{0.2cm}
\includegraphics[scale=0.4]{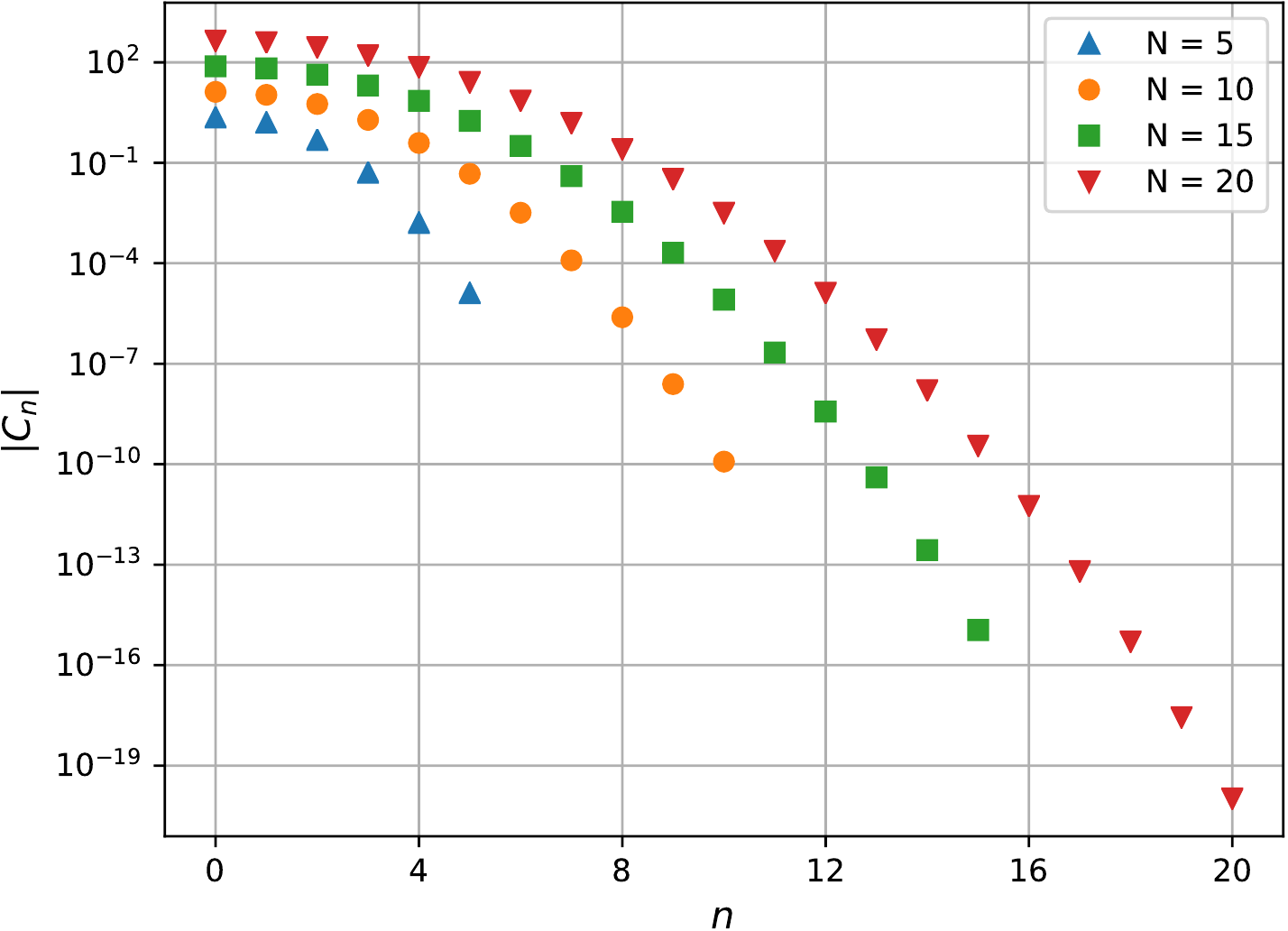}
\end{center}
\caption{The magnitudes~$|C_n|$ of the complex weights for $0\le n\le N$ 
and four choices of~$N$ using parabolic (left) and hyperbolic (right) contours.}
\label{fig: quad coeffs}
\end{figure}

After truncating the infinite sum as before to obtain a practical quadrature 
approximation~\eqref{eq: QhN}, the overall 
error~$\varepsilon_N=\varepsilon_N(\mu,\phi,h,\delta)$ is of order
\[
M(r_\delta;z)\exp(-2\pi r_\delta/h)+M(-s_\delta;z)\exp(-2\pi s_\delta/h)
+\exp\bigl(\mu(1-\cosh(Nh)\sin\phi)\bigr).
\]
We again seek to minimise the overall quadrature error by balancing these three 
terms. Neglecting $\delta$ leads to the equations
\[
\frac{\pi(2\phi-\pi)}{h}=\mu-\frac{2\pi\phi}{h} 
    =\mu\bigl(1-\cosh(Nh)\,\sin\phi\bigr), 
\]
which imply that
\[
\mu=\frac{\pi(4\phi-\pi)}{h}\quad\text{and}\quad
\cosh(Nh)=\frac{2\phi}{(4\phi-\pi)\sin\phi}.
\]
In this way,
\begin{equation}\label{eq: hyperbolic mu h}
\mu=\frac{\pi(4\phi-\pi)}{a(\phi)}\,N
\quad\text{and}\quad h=\frac{a(\phi)}{N}
\end{equation}
where
\[
a(\phi)=\arcosh\biggl(\frac{2\phi}{(4\phi-\pi)\sin\phi}\biggr)
\quad\text{for $\frac{\pi}{4}<\phi<\frac{\pi}{2}$,}
\]
and with these choices of $\mu=\mu_N(\phi)$ and $h=h_N(\phi)$ the quadrature 
error satisfies
\[
\epsilon_N(\phi)=O\bigl(\exp(-b(\phi)N)\bigr)
\quad\text{with}\quad b(\phi)=\frac{\pi(\pi-2\phi)}{a(\phi)}.
\]
The function~$b(\phi)$ has a unique maximum value for~$\phi\in[\pi/4,\pi/2]$
when $\phi=\phi_\star\doteq1{\cdot}17210$, so $\varepsilon_N(\phi_\star)$ is 
of order $\exp\bigl(-b(\phi_*)N\bigr)\approx10.13^{-N}$, indicating somewhat 
faster convergence than for the parabolic contour.  Taking 
$\phi=\phi_\star$ in~\eqref{eq: hyperbolic mu h} gives the optimized parameter 
values $\mu_\star\doteq 4{\cdot}49198\times N$~and 
$h_\star\doteq 1{\cdot}08180/N$.

Recalling the form of~$g$ from~\eqref{eq: f g}, and noting that
\[
\frac{h_\star w'(u)}{2\pi i}=\frac{\mu_\star h_\star}{2\pi}\,\cos(iu-\phi_\star)
    =\frac{4\phi_\star-\pi}{2}\,\cos(iu-\phi_\star),
\]
we obtain once again a sum of the form~\eqref{eq: QN opt} but now with
\begin{equation}\label{eq: A Cn hyperbola}
A=2\phi_\star-\frac{\pi}{2}\quad\text{and}\quad 
C_n=e^{w(nh_\star)}\cos(inh_\star-\phi_\star).
\end{equation}
The coefficients again satisfy~\eqref{eq: conj} so it suffices to store 
$w(nh_\star)$~and $C_n$ for~$0\le n\le N$. 

\Cref{fig: quad points} compares the locations of the quadrature 
points~$w(nh_\star)$ on the parabolic and hyperbolic contours in the complex 
plane, for four choices of~$N$.  The magnitudes~$|C_n|$ of the correponding 
coefficients are compared in \cref{fig: quad coeffs}.  

\begin{remark}\label{remark: rational}
Recalling the form of~$f(w;z)$ from~\eqref{eq: f w z}, we see that 
\eqref{eq: QN opt} provides a rational approximation to the contour integral,
\[
Q_{\star,N}(f;z)=\sum_{n=-N}^N\frac{R_n}{z-P_n},
\]
with poles~$P_n=w(nh_\star)^\alpha$ and associated residues 
$R_n=-AC_nw(nh_\star)^{\alpha-\beta}$.  Of course, these poles all lie outside 
the sector $\alpha\pi<|\arg z|\le\pi$ where $Q_{\star,N}(f;z)$ is used.
\end{remark}

\subsection{Extending the method to handle a pole in the integrand}

Still assuming $0<\alpha<1$, we suppose now that $|\theta|\le\alpha\pi$ and
recall that $w^\alpha-z=0$ when $w=\gamma$; see~\eqref{eq: gamma}.  Write
\[
\frac{w^{\alpha-\beta}}{w^\alpha-z}=\frac{w^{\alpha-\beta}\rho(w;z)}{w-\gamma}
\quad\text{where}\quad\rho(w;z)=\frac{w-\gamma}{w^\alpha-z},
\]
and observe that $\rho(w;z)$ has a removable singularity at~$w=\gamma$ 
with~$\rho(\gamma;z)=\gamma^{1-\alpha}/\alpha$. We define
\begin{equation}\label{eq: f1 w z}
f_1(w;z)=\frac{w^{\alpha-\beta}\rho(w;z)
-\gamma^{\alpha-\beta}\rho(\gamma;z)}%
{w-\gamma}=\frac{w^{\alpha-\beta}}{w^\alpha-z}
-\frac{\alpha^{-1}\gamma^{1-\beta}}{w-\gamma},
\end{equation}
so that
\[
f(w;z)=\frac{w^{\alpha-\beta}}{w^\alpha-z}
    =\frac{\alpha^{-1}\gamma^{1-\beta}}{w-\gamma}+f_1(w;z)
\]
and hence by~\eqref{eq: f w z},
\begin{equation}\label{eq: E gamma f1}
E_{\alpha,\beta}(z)=\alpha^{-1}\gamma^{1-\beta}\exp(\gamma)
    +\frac{1}{2\pi i}\int_{\mathcal{C}}e^wf_1(w;z)\,dw.
\end{equation}
Since $f_1$ is analytic in the whole cut plane, we can deform~$\mathcal{C}$ into 
a parabolic or hyperbolic contour of the type considered above, and use the 
corresponding quadrature sum~$Q_{\star,N}(f_1;z)$ to approximate the integral.  
Note that $f_1(w;z)$ is not a rational function of~$z$ since 
$\gamma=z^{1/\alpha}$.

To avoid roundoff problems evaluating $f_1(w;z)$ when~$w$ is close to~$\gamma$,
let $\epsilon=(w-\gamma)/\gamma$ so that $w=\gamma(1+\epsilon)$, and define
\begin{equation}\label{eq: psi1}
\psi_{1,\alpha}(\epsilon)=\frac{(1+\epsilon)^\alpha-1}{\epsilon}
    =\sum_{k=1}^\infty\binom{\alpha}{k}\epsilon^{k-1}
\end{equation}
and
\begin{equation}\label{eq: psi2}
\psi_{2,\alpha}(\epsilon)
    =\frac{(1+\epsilon)^\alpha-(1+\alpha\epsilon)}{\epsilon^2}
    =\sum_{k=2}^\infty\binom{\alpha}{k}\epsilon^{k-2}
\end{equation}
for $|\epsilon|<1$. The functions $\psi_{1,\alpha}$ and $\psi_{2,\alpha}$ may be 
evaluated to high accuracy when~$|\epsilon|$ is small by suitable 
truncation of these Taylor expansions or by employing the library functions 
\texttt{expm1}~and \texttt{log1p}. Since 
$w^\alpha-z=z\epsilon\psi_{1,\alpha}(\epsilon)
=z\epsilon[\alpha+\epsilon\psi_{2,\alpha}(\epsilon)]$ we find that
\[
f_1(w;z)=\frac{\psi_{1,\alpha-\beta}(\epsilon)
-\alpha^{-1}\psi_{2,\alpha}(\epsilon)}{\gamma^\beta\psi_{1,\alpha}(\epsilon)}
\quad\text{with}\quad
f_1(\gamma;z)=\frac{1+\alpha-2\beta}{2\alpha\gamma^\beta}.
\]

Summarizing: for~$0<\alpha<1$, the Mittag-Leffler 
function~$E_{\alpha,\beta}(z)$ is approximated by
\begin{equation}\label{eq: E alpha beta N}
E_{\alpha,\beta,N}(z)=\begin{cases}
    Q_{\star,N}(f;z)&\text{if $\alpha\pi<|\arg z|\le\pi$,}\\
   \alpha^{-1}z^{(1-\beta)/\alpha}\exp(z^{1/\alpha})+Q_{\star,N}(f_1;z)&
\text{if $|\arg z|\le\alpha\pi$.}
\end{cases}
\end{equation}
We will use superscripts ``par'' and ``hyp'' to distinguish between the 
parabolic and hyperbolic cases where needed, so that
\begin{equation}\label{eq: error behaviour}
E_{\alpha,\beta}(z)=E^{\mathrm{par}}_{\alpha,\beta,N}(z)+O(8.12^{-N})
\quad\text{and}\quad
E_{\alpha,\beta}(z)=E^{\mathrm{hyp}}_{\alpha,\beta,N}(z)+O(10.13^{-N}).
\end{equation}

\begin{figure}
\begin{center}
\includegraphics[scale=0.44]{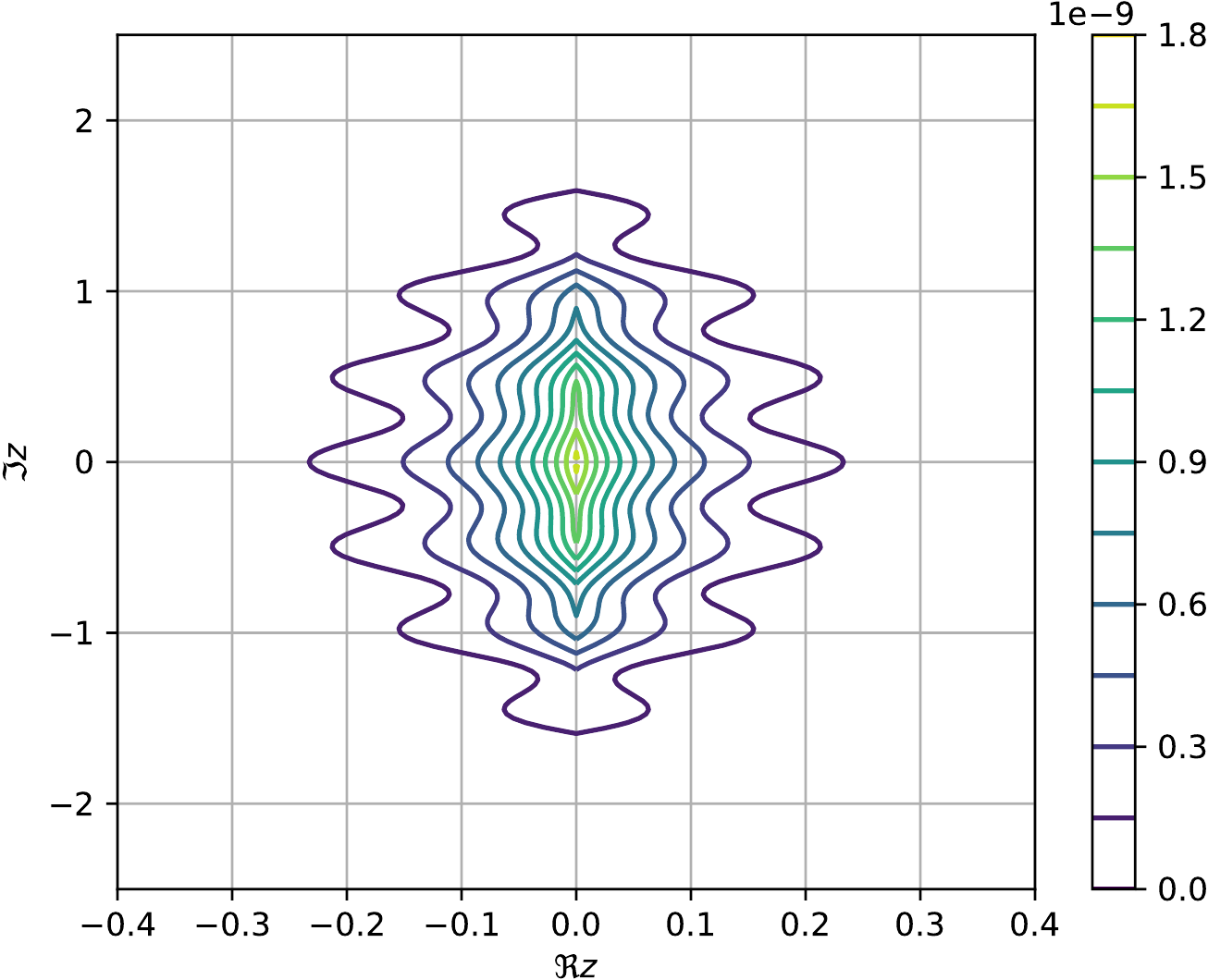} 
\hspace{0.2cm}
\includegraphics[scale=0.44]{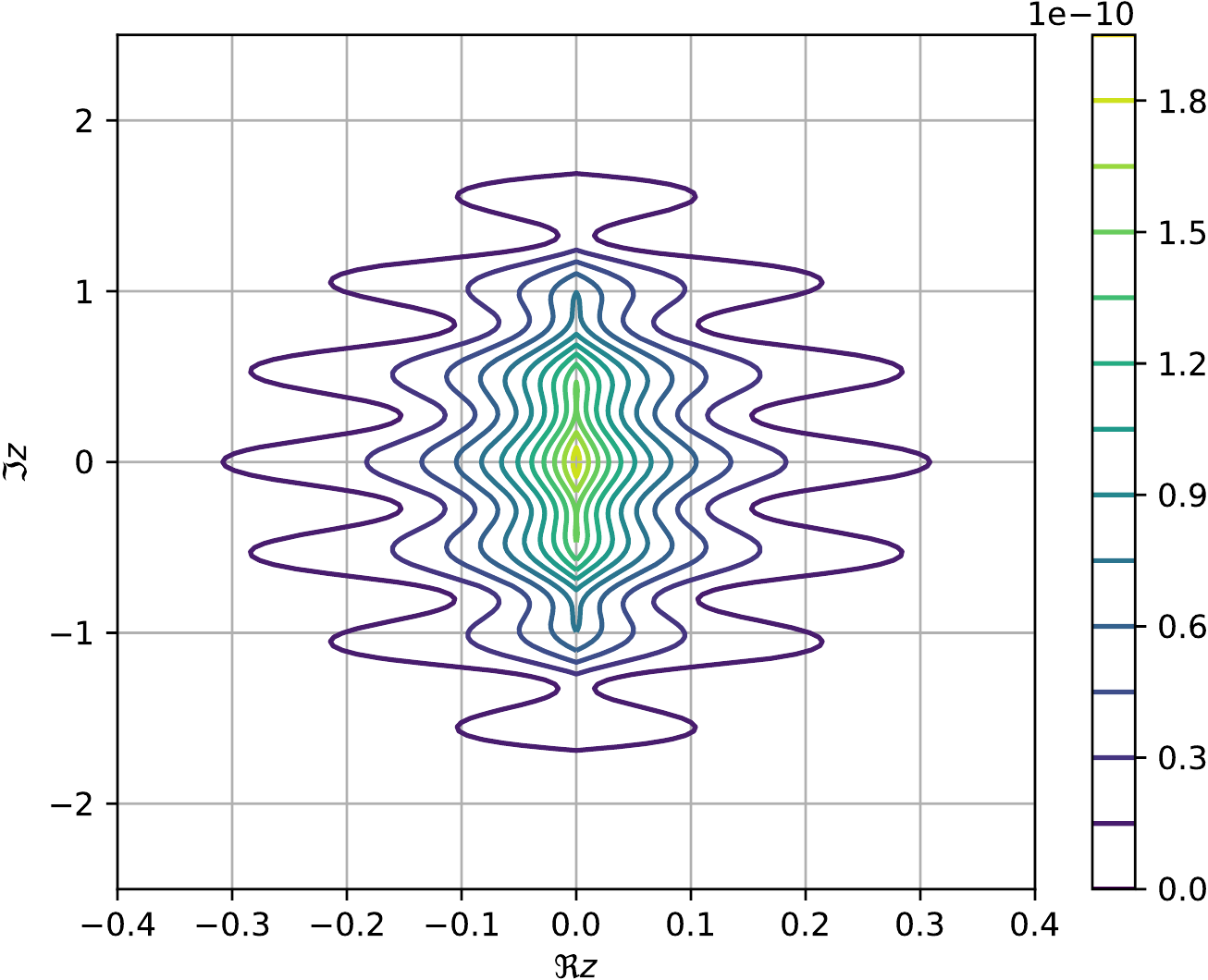} 
\end{center}
\caption{Magnitude of the error computing $E_{1/2}(z)$ for a quadrature 
rule with~$N=10$ using parabolic (left) and hyperbolic (right) contours.}
\label{fig: quad error}
\end{figure}

\begin{example}
We used both types of contour to compute $E_{1/2}(z)$ with~$N=10$ and obtained 
for the errors the patterns shown in \cref{fig: quad error}.  As expected, the 
hyperbolic contour gave somewhat more accurate results, with the error 
everywhere smaller than $1{\cdot}9\times10^{-10}$, compared 
to~$1{\cdot}7\times10^{-9}$ for the parabolic contour.  In both cases, the error 
is concentrated near the origin.  At the origin itself, the computed value
of~$E_\alpha(0+i0)$ is $\mathtt{NaN}$ due to the 
factor~$z^{(1-\beta)/\alpha}=(0+i0)^{(0+i0)}$ in~\eqref{eq: E alpha beta N}.  
\end{example}

\begin{remark}
When~$z$ is close to zero, it is usually best to compute~$E_{\alpha,\beta}(z)$ 
by truncating its Taylor expansion~\eqref{eq: gen E defn}.  In particular,
for~$\beta>1$ the factor~$z^{(1-\beta)/\alpha}$ in~\eqref{eq: E alpha beta N} 
blows up as~$z\to0$.
\end{remark}

\section{Evaluation on the real line}\label{sec: real line}

For $z$ on the real line, the functions from \eqref{eq: f w z}~and 
\eqref{eq: f1 w z} satisfy
\[
f(\bar w;x)=\overline{f(w;x)}
\quad\text{and}\quad
f_1(\bar w;x)=\overline{f_1(w;x)}
\quad\text{for $u$, $x\in\mathbb{R}$,}
\]
so it follows using the properties~\eqref{eq: conj} of the quadrature points 
and coefficients that we can nearly halve the number of function evaluations
required to approximate the contour integral because
\[
Q_{\star,N}(f;x)=A\biggl(C_0f\bigl(w(0);x\bigr)+2\sum_{n=1}^N
    \Re\bigl[C_nf\bigl(w(nh_\star);x\bigr)\bigr]\biggr).
\]

Assume $0<\alpha<1$ and $x>0$.  Taking $z=-x$ we have $\theta=\pi$ so 
\begin{equation}\label{eq: E -x}
E_{\alpha,\beta,N}(-x)=Q_{\star,N}(f;-x).
\end{equation}
However, taking $z=x$ means that $\theta=0$ and $\gamma=x^{1/\alpha}$, so
\begin{equation}\label{eq: E +x}
E_{\alpha,\beta,N}(x)=\alpha^{-1}x^{(1-\beta)/\alpha}\exp(x^{1/\alpha})
    +Q_{\star,N}(f_1;x).
\end{equation}
In fact, \eqref{eq: E +x} is valid for $0<\alpha<2$ because $w=x^{1/\alpha}$ is 
still the only solution to~$w^\alpha-x=0$ with $|\arg w|\le\pi$.

\begin{figure}
\begin{center}
\includegraphics[scale=0.42]{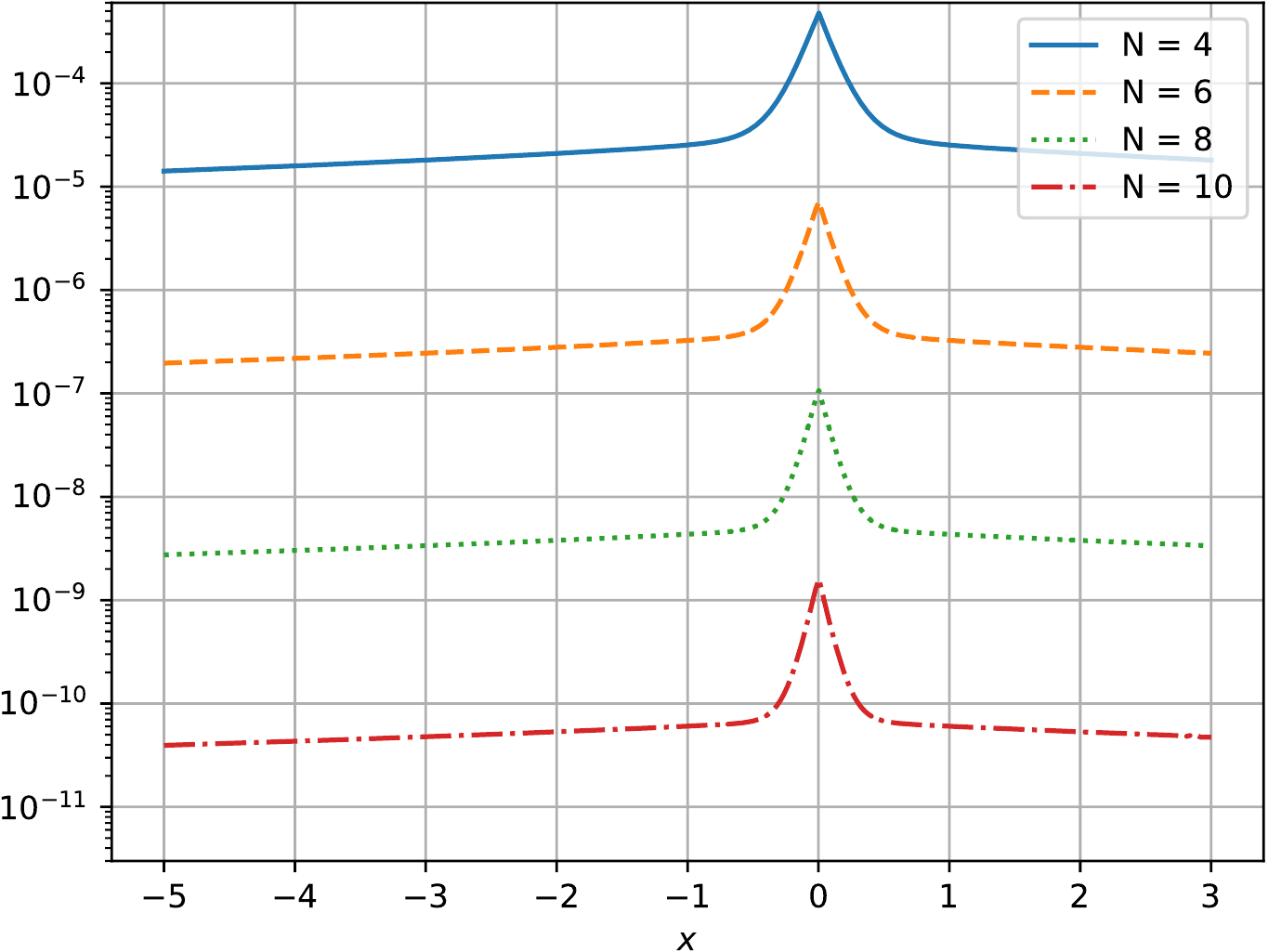}
\hspace{0.2cm}
\includegraphics[scale=0.42]{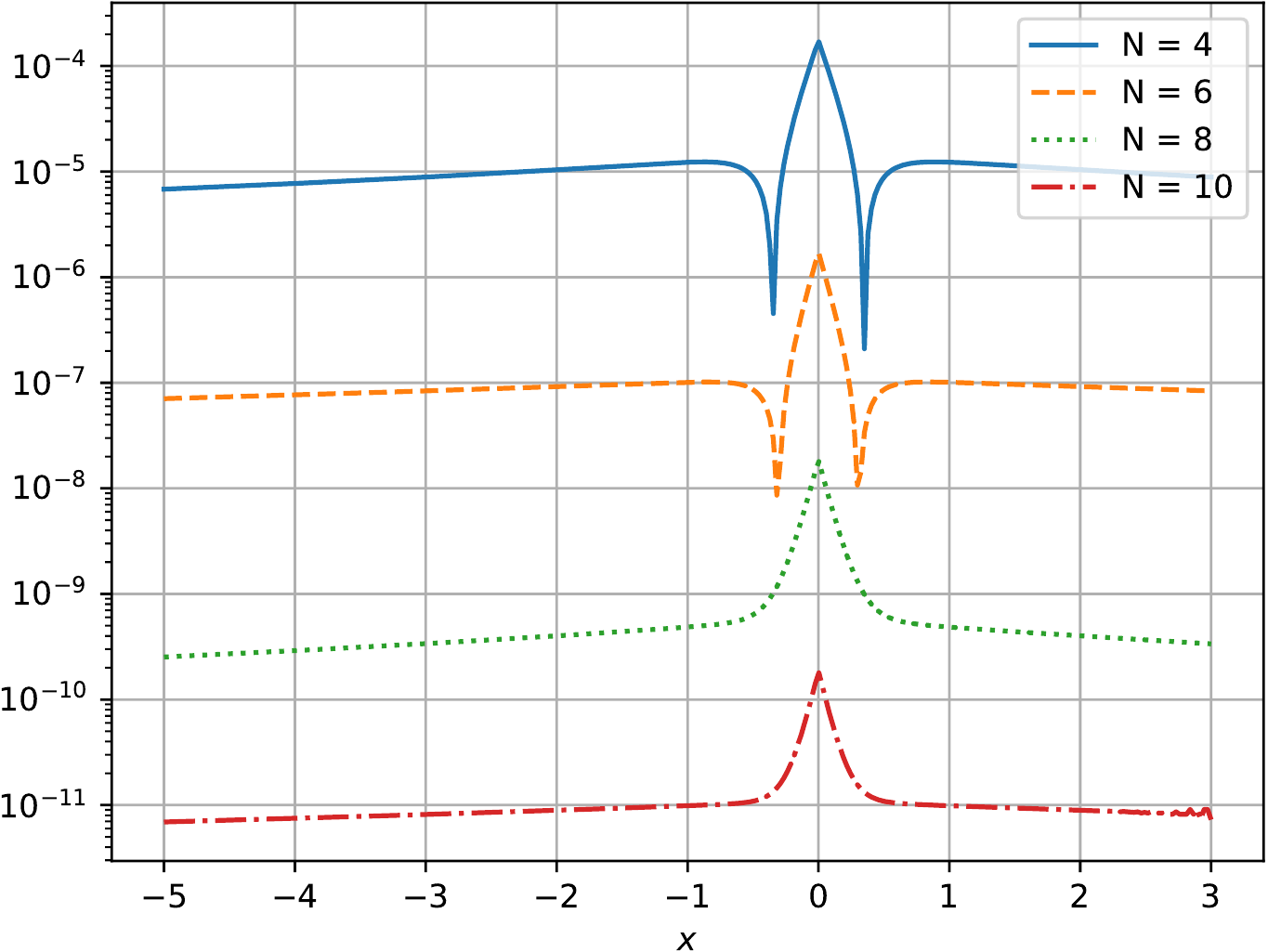}
\end{center}
\caption{The absolute error~$|E_{1/2,1}(-x)-E_{1/2,1,N}(-x)|$ for increasing 
choices of~$N$, using parabolic (left) and hyperbolic (right) 
contours.}\label{fig: quad conv}
\end{figure}

\begin{example}
\Cref{fig: quad conv} shows the absolute errors in~$E_{1/2,1,N}(x)$ 
for~$-5\le x\le3$ on a log scale.  The convergence behaviour is consistent
with the predictions~\eqref{eq: error behaviour}
from the analysis above.  Similar calculations reveal that the accuracy  
does not vary much with respect to~$\alpha\in(0,1]$.  As observed earlier in 
\cref{fig: quad error} the maximum error occurs at~$x=0$,  and since 
$f(w;0)=w^{-\beta}$ and $E_{\alpha,\beta}(0)=1/\Gamma(\beta)$ are 
independent of~$\alpha$, the difference $Q_{\star,N}(f;0)-1/\Gamma(\beta)$ 
provides a readily computable measure of the expected accuracy 
of~$E_{\alpha,\beta,N}(-x)$ for $x>0$~and $0<\alpha\le1$.  On 
the positive half-line, we find that as $x$ increases from~$0$ the computed 
value~$E_{\alpha,\beta,N}(x)$ is quickly dominated by the first term 
in~\eqref{eq: E +x}.  For instance,
$E_{1/2}(3)\approx2\exp(3^2)\doteq16,206$.
\end{example}

Suppose now that $1<\alpha<2$, $z=-x$~and $x>0$. The equation 
$w^\alpha+x=0$ then has two solutions 
$w=\gamma_\pm=x^{1/\alpha}\exp(\pm i\pi/\alpha)$ satisfying
\[
-\pi<\arg\gamma_-<-\frac{\pi}{2}\quad\text{and}\quad
\frac{\pi}{2}<\arg\gamma_+<\pi.
\]
As an alternative to applying the identity~\eqref{eq: E alpha/m}
we can deal with the singularities~$\gamma_\pm$ using a similar, albeit more 
complicated, approach to the one that led to~\eqref{eq: E +x}.  Write
\[
\frac{1}{w^\alpha+x}=\varphi(w;x)\biggl(
    \frac{1}{w-\gamma_+}+\frac{1}{w-\gamma_-}\biggr),
\]
where the function
\[
\varphi(w;x)=\frac{(w-\gamma_+)(w-\gamma_-)}%
{(w^\alpha+x)(2w-\gamma_+-\gamma_-)}
\]
has removable singularities at~$w=\gamma_\pm$ 
with~$\varphi(\gamma_\pm;x)=\alpha^{-1}\gamma_\pm^{1-\alpha}$. Define 
\[
f_\pm(w;x)=\frac{w^{\alpha-\beta}\varphi(w;x)
-\gamma_\pm^{\alpha-\beta}\varphi(\gamma_\pm;x)}{w-\gamma_\pm}
=\frac{w^{\alpha-\beta}(w-\gamma_\mp)}{(w^\alpha+x)(2w-\gamma_+-\gamma_-)}
    -\frac{\alpha^{-1}\gamma_\pm^{1-\beta}}{w-\gamma_\pm}
\]
so that
\[
w^{\alpha-\beta}\,\frac{\varphi(w;x)}{w-\gamma_\pm}
    =\frac{\alpha^{-1}\gamma_\pm^{1-\beta}}{w-\gamma_\pm}+f_\pm(w;x).
\]
In this way, letting $f_2(w;x)=f_+(w;x)+f_-(w;x)$, it follows that
\[
f(w;-x)=\frac{w^{\alpha-\beta}}{w^\alpha+x}
    =\frac{1}{\alpha}\biggl(\frac{\gamma_+^{1-\beta}}{w-\gamma_+}
    +\frac{\gamma_-^{1-\beta}}{w-\gamma_-}\biggr)+f_2(w;x).
\] 
Putting $z=-x$ in the integral representation~\eqref{eq: Wiman integral},
we see that
\begin{equation}\label{eq: E neg 1<alpha<2}
E_{\alpha,\beta}(-x)=\frac{1}{\alpha}\Bigl(\gamma_+^{1-\beta}e^{\gamma_+}
    +\gamma_-^{1-\beta}e^{\gamma_-}\Bigr)
    +\frac{1}{2\pi i}\int_{\mathcal{C}}e^w f_2(w;x)\,dw,
\end{equation}
with $f_2(w;x)$ analytic for all~$w$ in the cut plane and satisfying
$f_2(\bar w;x)=\overline{f_2(w;x)}$, and with
\begin{equation}\label{eq: cosine term}
\gamma_+^{1-\beta}e^{\gamma_+}+\gamma_-^{1-\beta}e^{\gamma_-}
    =2x^{(1-\beta)/\alpha}\exp\biggl(x^{1/\alpha}\cos\frac{\pi}{\alpha}\biggr)\,
    \cos\biggl((1-\beta)\frac{\pi}{\alpha}
    +x^{1/\alpha}\sin\frac{\pi}{\alpha}\biggr).
\end{equation}

To compute $f_\pm(w;x)$ for~$w$ close to~$\gamma_\pm$, we let 
$\epsilon_\pm=\frac{w-\gamma_\pm}{\gamma_\pm}$ so that
\[
w=\gamma_\pm(1+\epsilon_\pm),\quad
w-\gamma_\pm=\gamma_\pm\epsilon_\pm,\quad
2w-\gamma_\pm-\gamma_\mp=w-\gamma_\mp+\epsilon_\pm\gamma_\pm,
\]
and
\[
w^\alpha+x=\gamma_\pm^\alpha\epsilon_\pm\psi_{1,\alpha}(\epsilon_\pm)
    =\gamma_\pm^\alpha\epsilon_\pm
    \bigl[\alpha+\epsilon_\pm\psi_{2,\alpha}(\epsilon_\pm)\bigr].
\]
We find that
\[
f_\pm(w;x)=\frac{(w-\gamma_\mp)[\psi_{1,\alpha-\beta}(\epsilon_\pm)
-\alpha^{-1}\psi_{2,\alpha}(\epsilon_\pm)]
-\gamma_\pm\alpha^{-1}\psi_{1,\alpha}(\epsilon_\pm)}%
{\gamma_\pm^\beta\psi_{1,\alpha}(\epsilon_\pm)
(w-\gamma_\mp+\epsilon_\pm\gamma_\pm)},
\]
and in particular,
\[
f_\pm(\gamma_\pm;x)=\frac{
(1+\alpha-2\beta)(\gamma_\pm-\gamma_\mp)-2\gamma_\pm}%
{2\alpha\gamma_\pm^\beta(\gamma_\pm-\gamma_\mp)}.
\]
Similarly,
\[
f_\mp(w;x)=\frac{\gamma_\pm^{1-\beta}(1+\epsilon_\pm)^{\alpha-\beta}}%
{\psi_{1,\alpha}(\epsilon)(w-\gamma_\mp+\epsilon_\pm\gamma_\pm)}
-\frac{\alpha^{-1}\gamma_\mp^{1-\beta}}{w-\gamma_\mp},
\]
and in particular,
\[
f_\mp(\gamma_\pm;x)=\frac{\gamma_\pm^{1-\beta}-\gamma_\mp^{1-\beta}}%
{\alpha(\gamma_\pm-\gamma_\mp)}=x^{-\beta/\alpha}\,
\frac{\sin\pi(1-\beta)/\alpha}{\alpha\sin\pi/\alpha}.
\]
\section{Approximations of Pad\'e type}\label{sec: Pade}
For $0<\alpha\le1$, consider the problem of approximating 
$E_{\alpha,\beta}(-x)$ by a rational function $p(x)/q(x)$ where
\[
p(x)=\sum_{j=0}^rp_jx^j\quad\text{and}\quad q(x)=\sum_{j=0}^rq_jx^j.
\]
Zeng and Chen~\cite{ZengChen2015} proposed a two-point Pad\'e scheme, 
with the coefficients determined by the conditions 
\begin{equation}\label{eq: Pade condition}
E_{\alpha,\beta}(-x)=\frac{p(x)}{q(x)}+\begin{cases}
    O(x^m)&\text{as $x\to0$,}\\
    O(x^{-n})&\text{as $x\to\infty$,}
\end{cases}
\end{equation}
with appropriate choices of $m$~and $n$.  They obtained closed-form expressions 
for $p_j$~and $q_j$ in terms of $\alpha$~and $\beta$ for some small choices 
of~$r$. Let
\[
a(x)=\sum_{k=0}^{m-1}a_kx^k\quad\text{and}\quad
b(x)=\sum_{k=1}^{n-1}b_kx^k,
\]
where, recalling \eqref{eq: Gamma reflection},
\[
a_k=\smash[t]{\frac{(-1)^k}{\Gamma(\beta+k\alpha)}}
\quad\text{and}\quad
b_k=\frac{(-1)^{k-1}}{\Gamma(\beta-k\alpha)}
    =\frac{(-1)^k}{\pi}\,\sin\pi(k\alpha-\beta)\,\Gamma(1+k\alpha-\beta).
\]
With the above definitions, \eqref{eq: E defn}~and 
\eqref{eq: asymptotics} imply that
\[
E_{\alpha,\beta}(-x)=\begin{cases}
    a(x)+O(x^m)&\text{as $x\to0$,}\\
    b(x^{-1})+O(x^{-n})&\text{as $x\to\infty$,}
\end{cases}
\]
so the conditions in~\eqref{eq: Pade condition} hold iff
\[
\frac{p(x)}{q(x)}-a(x)=O(x^m)\quad\text{as $x\to0$,}
\]
and
\[
\frac{x^{-r}p(x)}{x^{-r}q(x)}-b(x^{-1})=O(x^{-n})\quad\text{as $x\to\infty$,}
\]
or equivalently, assuming $q_0\ne0$,
\begin{equation}\label{eq: Pade x->0}
p(x)-a(x)q(x)=O(x^m)\quad\text{as $x\to0$,}
\end{equation}
and, assuming $q_r\ne0$,
\begin{equation}\label{eq: Pade x->oo}
x^{-r}\bigl[p(x)-b(x^{-1})q(x)\bigr]=O(x^{-n})\quad\text{as $x\to\infty$.}
\end{equation}
Setting the coefficient of~$x^k$ on the left-hand side 
of~\eqref{eq: Pade x->0} to zero for~$0\le k\le m-1$, and the coefficient 
of~$x^{-k}$ on the left-hand side of~\eqref{eq: Pade x->oo} to zero 
for~$0\le k\le n-1$, provides us with $m+n$~equations for the $2(r+1)$ unknown 
coefficients of the polynomials $p$~and $q$.  However, since the ratio~$p/q$ is 
unchanged if the numerator and denominator are multiplied by a common non-zero 
factor, we are free to impose a scaling condition that reduces the number of 
degrees of freedom to~$2r+1$, suggesting that we should require
\[
m+n=2r+1.
\]
Furthermore, since $b(0)=0$, the left-hand side of~\eqref{eq: Pade x->oo} 
equals $p_r+O(x^{-1})$ so we must have $p_r=0$.

If $m\ge r+1$, then $n=2r+1-m\le r$ and we find that
\[
\begin{aligned}
p_k-\sum_{j=0}^k a_{k-j}q_j&=0&&\text{for $0\le k\le r-1$,}\\
-\sum_{j=0}^r a_{k-j}q_j&=0&&\text{for $r\le k\le m-1$,}\\
p_k-\sum_{j=k+1}^rb_{j-k}q_j&=0&&\text{for $r-n+1\le k\le r-1$,}
\end{aligned}
\]
whereas if $m\le r$, then $n=2r+1-m\ge r+1$ and
\[
\begin{aligned}
p_k-\sum_{j=0}^k a_{k-j}q_j&=0&&\text{for $0\le k\le m-1$,}\\
-\sum_{j=0}^r b_{j-k}q_j&=0&&\text{for $-(r-m)\le k\le-1$,}\\
p_k-\sum_{j=k+1}^r b_{j-k}q_j&=0&&\text{for $0\le k\le r-1$.} 
\end{aligned}
\]
In either case, we obtain a $(2r)\times(2r+1)$ homogeneous linear 
system of the form
\[
C\mathbf{x}=\boldsymbol{0}\quad\text{where}\quad
\mathbf{x}=[p_0,p_1,\ldots,p_{r-1},q_0,q_1,\ldots,q_r]^\top.
\]
If the matrix~$C$ has full rank, then $\boldsymbol{x}$ must belong to a 
one-dimensional subspace of~$\mathbb{R}^{2r+1}$.  

Let $\mathbf{c}_j$ denote the $j$th column of~$C$.
Following Zeng and Chen~\cite{ZengChen2015} we can fix a solution by 
putting~$x_{r+1}=q_0=1$ to obtain a $(2r)\times(2r)$ linear system
\[
\tilde C\tilde{\mathbf{x}}=-\mathbf{c}_{r+1}
\]
where $\tilde 
C=[\mathbf{c}_1,\ldots,\mathbf{c}_r, \mathbf{c}_{r+2},
\ldots,\mathbf{c}_{2r+1}]$ and
$\tilde{\mathbf{x}}=[p_0,\ldots,p_{r-1},q_1,\ldots,q_r]^\top$.  
This linear system is badly conditioned, the more so the smaller 
the value of~$\alpha$ and the larger the value of~$r$, but in practice we 
observe an autocorrection phenomenon \cite{Litvinov2003,Luke1980} so that
the computed values of~$p(x)/q(x)$ can nevertheless provide an accurate 
approximation to~$E_{\alpha,\beta}(-x)$.

In fact, computing the singular-value decomposition (SVD),
\[
\tilde C=\tilde U\tilde S\tilde V^\top,\quad 
\tilde S=\diag(\tilde\sigma_1,\tilde\sigma_2,\ldots,\tilde\sigma_{2r}),\quad
\tilde\sigma_1\ge\tilde\sigma_2\ge\cdots\ge\tilde\sigma_{2r},
\]
so that $\tilde S(\tilde V^\top\tilde{\mathbf{x}})
=-\tilde U^\top\mathbf{c}_{r+1}$, we find that the last few singular values 
become very small but so do the corresponding components 
of~$\tilde U^\top\mathbf{c}_{r+1}$. Thus, the computed coefficients in the
solution vector~$\tilde{\mathbf{x}}$ are of moderate size.

Instead of fixing $q_0=1$, we can compute the SVD of the whole 
$(2r)\times(2r+1)$ matrix, $C=USV^T$.  Here, the diagonal matrix~$S$ is 
$(2r)\times(2r+1)$, and hence the orthogonal matrices $U$~and $V$ are 
$(2r)\times(2r)$~and $(2r+1)\times(2r+1)$, respectively.  Let 
$\mathbf{e}_j$ denote the $j$th standard basis vector 
in~$\mathbb{R}^{2r+1}$ and put 
\[
\mathbf{x}=V\mathbf{e}_{2r+1}=\text{last column of $V$,}
\]
so that $V^\top\mathbf{x}=\mathbf{e}_{2r+1}$. It follows that
$C\mathbf{x}=USV^T\mathbf{x}=US\mathbf{e}_{2r+1}=\mathbf{0}$ 
because the last column of~$S$ is zero.  Since $V$ is orthogonal, this choice
of~$\mathbf{x}$ results in the normalization
\[
\sum_{j=0}^{r-1}p_j^2+\sum_{j=0}^rq_j^2=1.
\]
Having computed $p(x)/q(x)$ in this way, if desired we can of course rescale 
the coefficients so that $q(0)=q_0=1$. 

A third option is to compute the LU factorization of~$C$ with partial pivoting,
$C=PLU$, and solve the $(2r)\times(2r+1)$, homogeneous, upper triangular
system $U\mathbf{x}=\mathbf{0}$ by setting $x_{2r+1}=q_r=1$ and using 
back substitution.  Provided $x_r=q_0\ne0$, the coefficients can again be 
rescaled so that $q_0=1$, if desired.  However, we found that $C$~and $U$ are 
nearly as badly conditioned as~$\tilde C$. For example, with $\alpha=0.2$, 
$\beta=1$, $m=9$~and $n=8$ the condition numbers of $\tilde C$, $C$~and $U$
are about $1{\cdot}39\times10^{13}$, $4{\cdot}83\times10^{12}$~and 
$1{\cdot}70\times10^{12}$, respectively, with the coefficients $p_j$~and 
$q_j$ in agreement to only about 4~significant figures.  Nevertheless, the 
values of $p(x)/q(x)$ computed using these three approaches agree to within 
about $6\times10^{-16}$ for~$0\le x<\infty$.

\begin{figure}
\begin{center}
\includegraphics[scale=0.6]{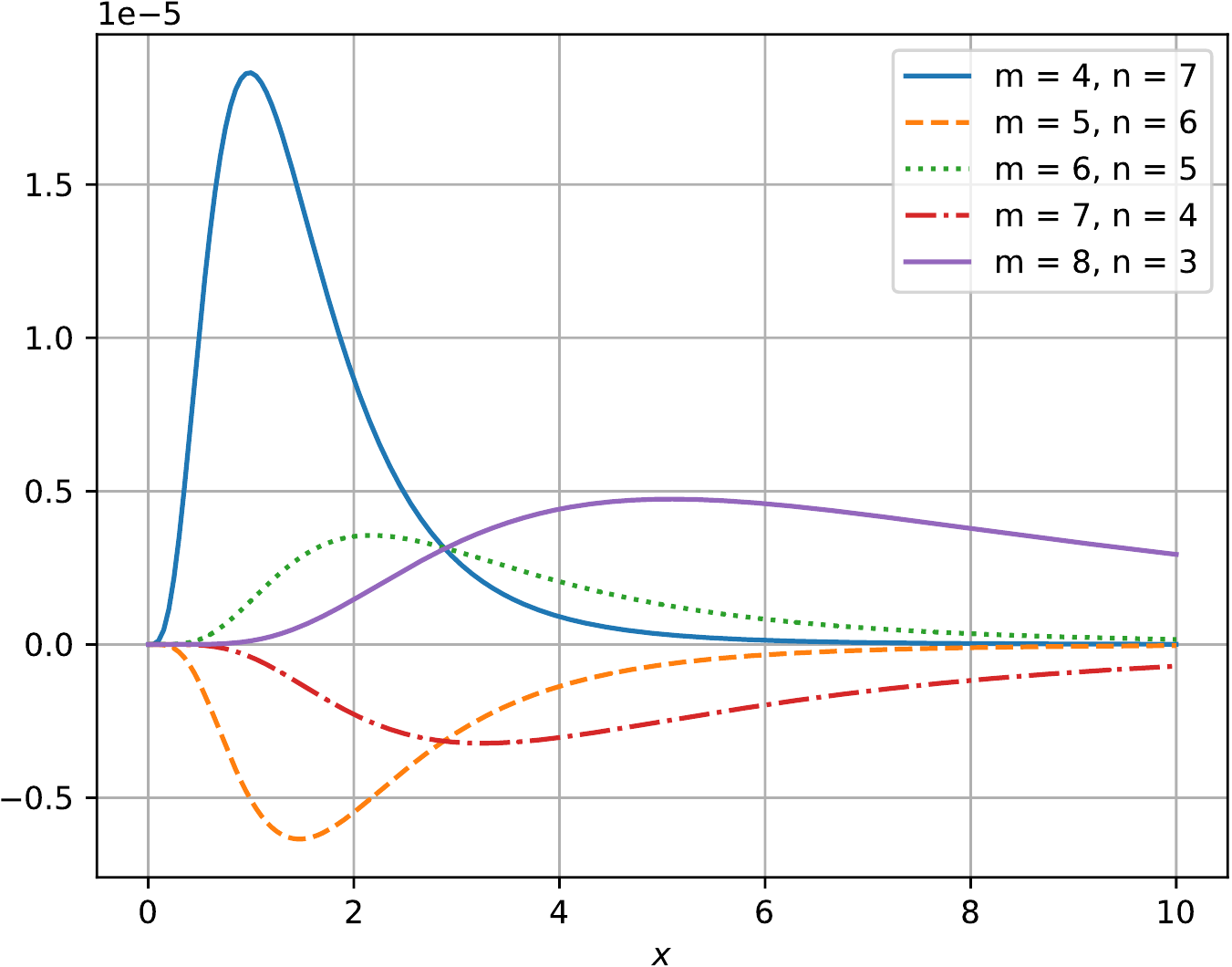} 
\end{center}
\caption{The error $p(x)/q(x)-E_{1/2}(-x)$ for different choices of $m$~and $n$
with $r=5$.}\label{fig: pade m n}
\end{figure}

\begin{figure}
\begin{center}
\includegraphics[scale=0.6]{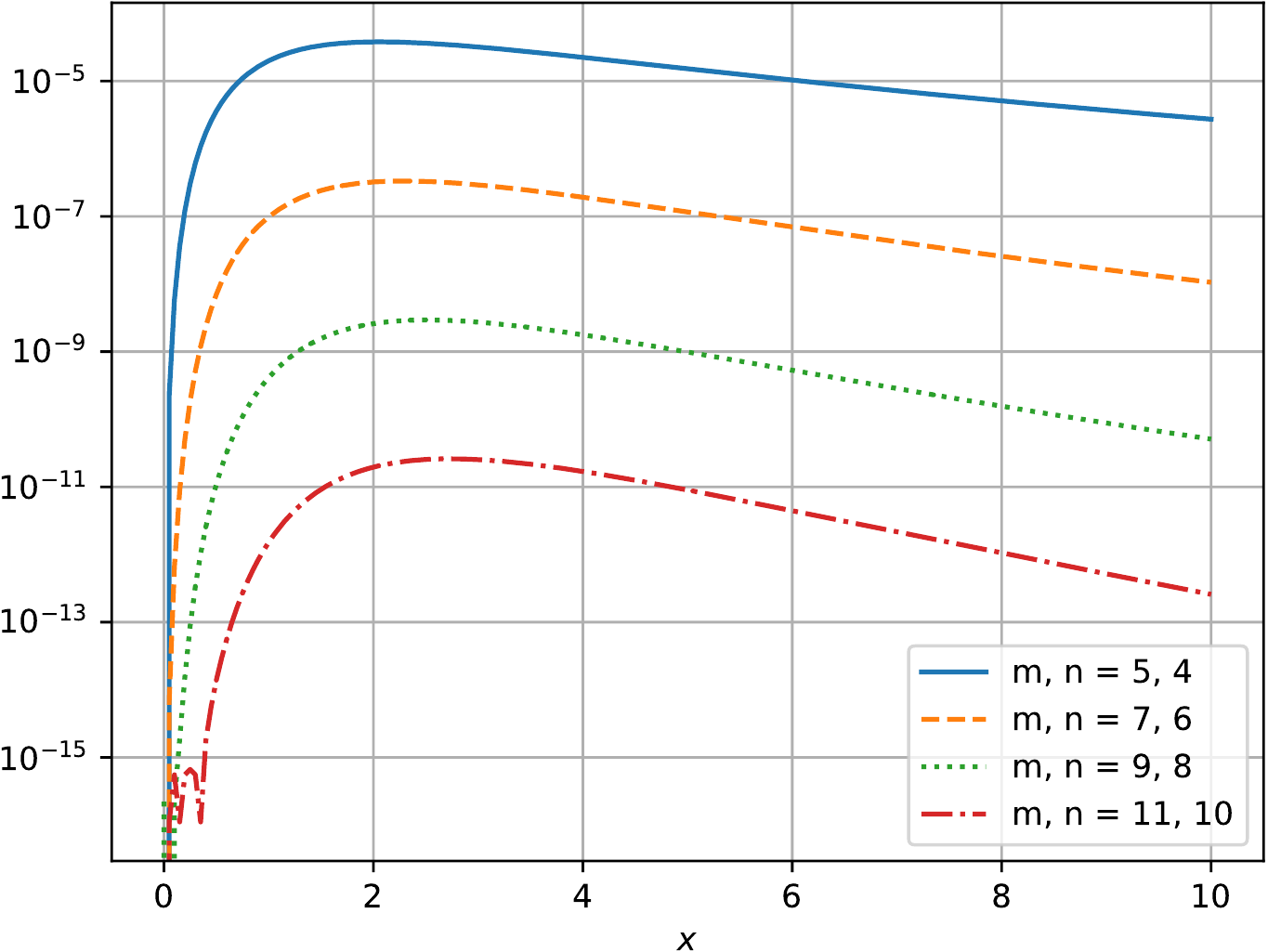} 
\end{center}
\caption{The absolute error $|E_{1/2}(-x)-p(x)/q(x)|$ when $m=r+1$ and $n=r$ 
for increasing choices of~$r$.}\label{fig: pade conv}
\end{figure}

\begin{figure}
\begin{center}
\includegraphics[scale=0.6]{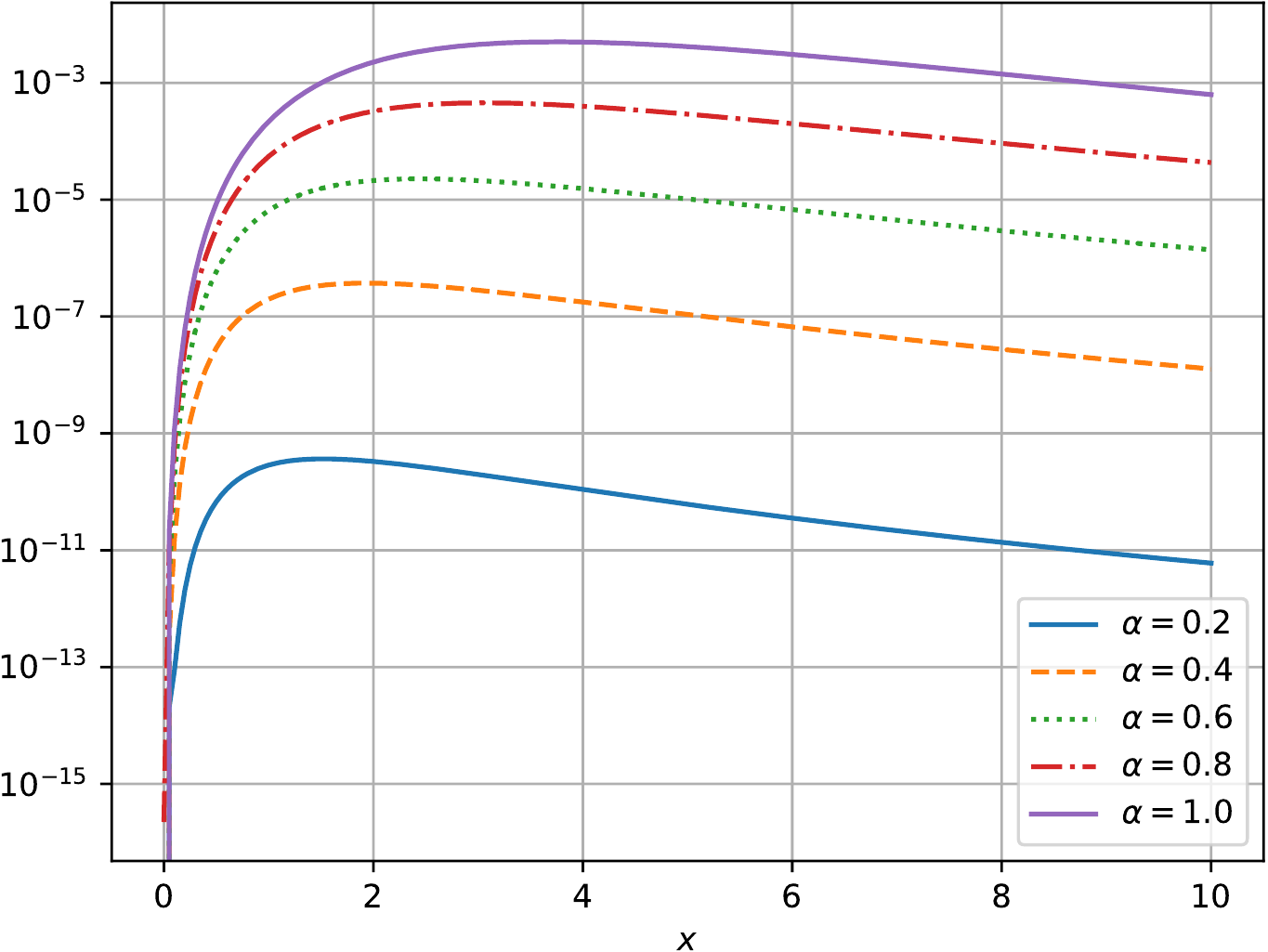} 
\end{center}
\caption{The absolute error $|p(x)/q(x)-E_\alpha(-x)|$ for different choices 
of~$\alpha$ with fixed polynomial degrees $m=6$~and $n=5$.}
\label{fig: pade alpha}
\end{figure}

\begin{example}
\Cref{fig: pade m n} shows the error in the Pad\'e approximation of 
$E_{1/2}(-x)$ for five choices of $m$~and $n$ such that $m+n=11$.  Not 
surprisingly, increasing $m$ improves the accuracy for small~$x$ at the expense 
of worse accuracy for large~$x$.  Increasing $n$ has the opposite effect,
and the accuracy for middling values of~$x$ is best when $m$~and $n$ are 
roughly equal. \Cref{fig: pade conv} plots the absolute error 
when $m=r+1$ and $n=r$, showing how the accuracy improves with increasing 
values of~$r$.
\end{example}

\begin{example}
\Cref{fig: pade alpha} illustrates how the accuracy of the Pad\'e approximation 
gets worse as~$\alpha$ increases while keeping $m$~and $n$ fixed; in this case,
$m=6$ and $n=5$.  
\end{example}

Let $\chi_1$, $\chi_2$, \ldots, $\chi_r$ denote the zeros of~$q$, so that
$q(z)=q_r\prod_{j=1}^r(z-\chi_j)$, and assume for simplicity that the roots are 
distinct. For~$z$ in some neighbourhood of the negative real axis,
consider the partial fraction expansion
\[
E_{\alpha,\beta}(z)\approx
\frac{p(-z)}{q(-z)}=\sum_{j=1}^r\frac{\varrho_j}{z+\chi_j}\quad\text{where}\quad
\varrho_j=-\frac{p(\chi_j)}{q'(\chi_j)}.
\]
We can compare this rational approximation with the one discussed in 
\cref{remark: rational}.

\begin{figure}
\begin{center}
\includegraphics[scale=0.44]{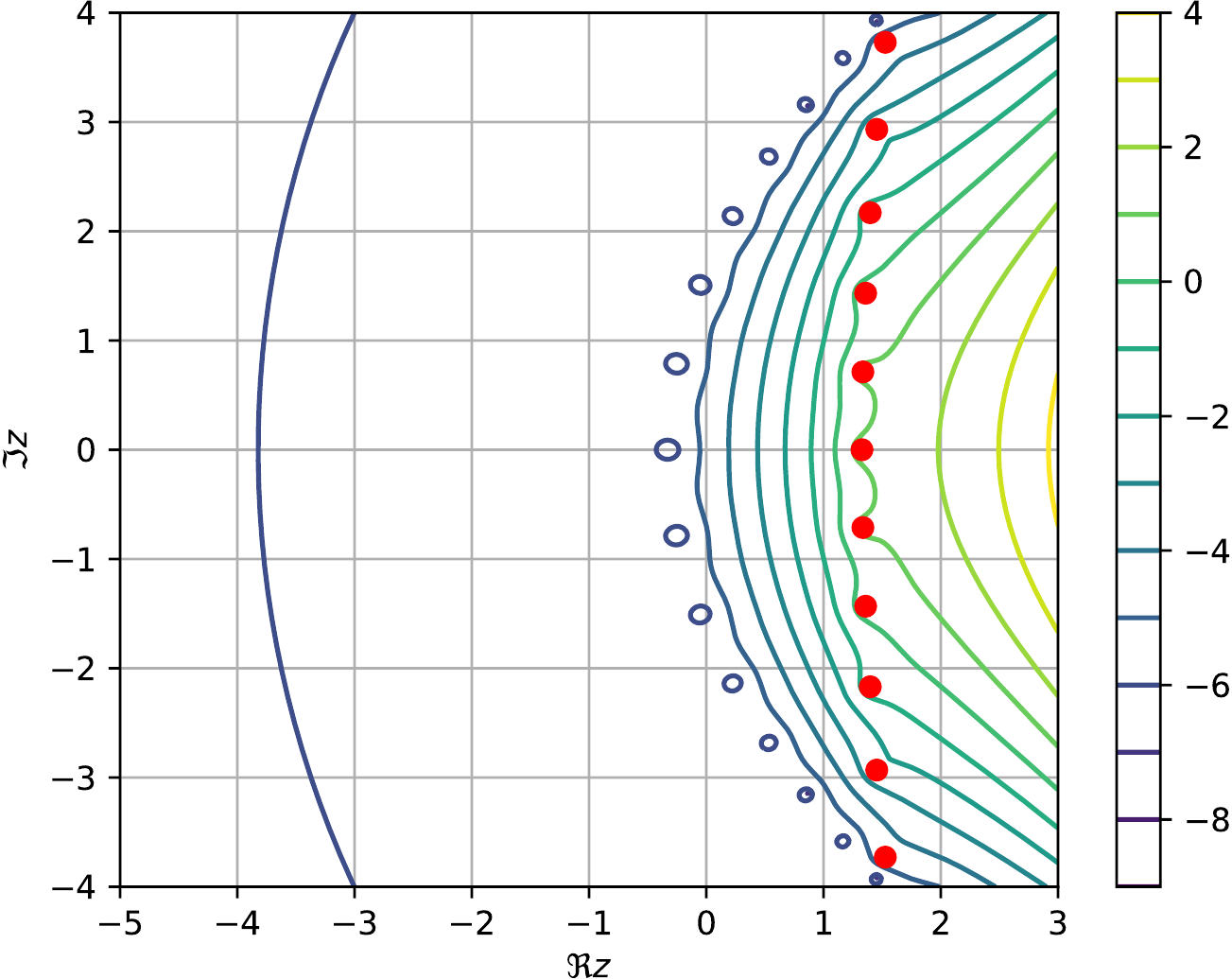}
\hspace{0.2cm}
\includegraphics[scale=0.44]{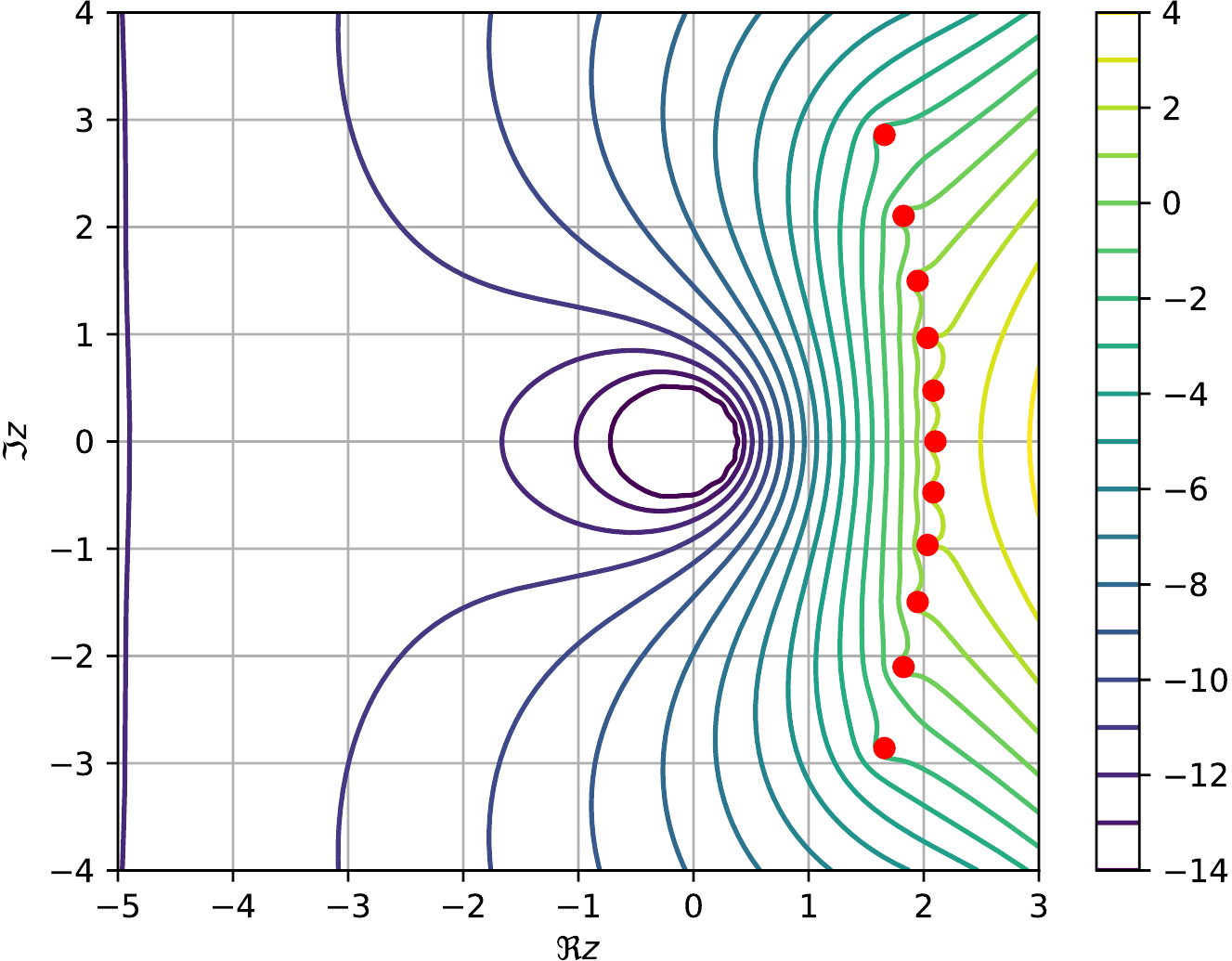}
\end{center}
\caption{Contour plots of 
$\log_{10}|E_{1/2}(z)-Q_{\star,N}^{\mathrm{hyp}}(f;z)|$ (left) for~$N=5$ and 
$\log_{10}|E_{1/2}(z)-p(-z)/q(-z)|$ (right) for $m=6$~and $n=5$.  The red dots 
mark the locations of the $11$~poles for each of the rational approximations.}
\label{fig: poles}
\end{figure}

\begin{example}
\Cref{fig: poles} shows contour plots of the base-10~logarithm of the absolute 
error for the two types of rational approximation to~$E_{1/2}(z)$ for
$-5\le\Re z\le 3$ and $-4\le\Im z\le4$, with the locations of the poles 
shown as red dots. The quadrature approximation using a hyperbola for~$N=5$ is 
compared with the Pad\'e approximation for $m=12$~and $n=11$, so that both sums 
have $r=2N+1=11$~terms.  Recall from~\eqref{eq: asymptotics} that 
$E_{1/2}(z)=2\exp(z^2)+O(z^{-1})$ as~$|z|\to\infty$ when~$|\arg z|\le\pi/2$, so
we cannot expect either approximation to work for~$\Re z>0$, except 
near~$z=0$.  However, both are effective for~$\Re z<0$.  In this case, the 
Pad\'e approximation~$p(-z)/q(-z)$ achieves higher accuracy than the quadrature 
approximation~$Q_{\star,N}^{\mathrm{hyp}}(z)$.
\end{example}
\section{Best approximation by rational functions}\label{sec: best approx}
For a real interval~$I$, let $\mathcal{R}^m_n(I)$ denote the set of rational 
functions~$p/q$ where $p:I\to\mathbb{R}$~and $q:I\to\mathbb{R}$ are polynomials
of degree at most $m$~and $n$, respectively. In light of the preceding results, 
it is natural to ask what is the best possible accuracy achievable when 
approximating $E_\alpha(-x)$ for~$0\le x<\infty$ by a function 
in~$\mathcal{R}^m_n\bigl([0,\infty)\bigr)$. For $\alpha=1$, that is, for the 
exponential function~$e^{-x}$, and for~$m=n$, this question was addressed by 
Nakatsukasa et al.~\cite[Section~6.8]{NakatsukasaEtAl2018} and we will adapt 
their approach.

\begin{figure}
\begin{center}
\includegraphics[scale=0.42]{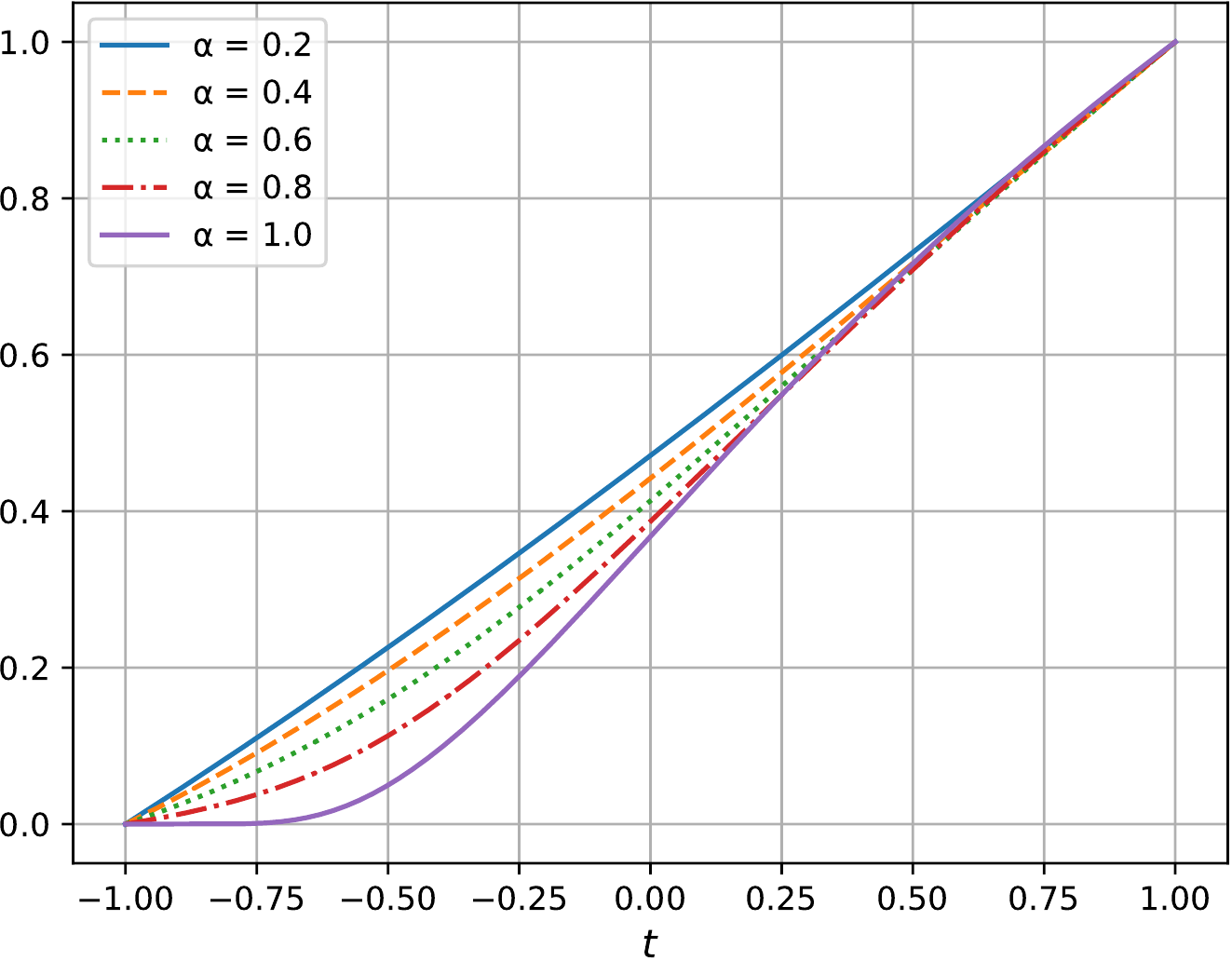} 
\hspace{0.2cm}
\includegraphics[scale=0.42]{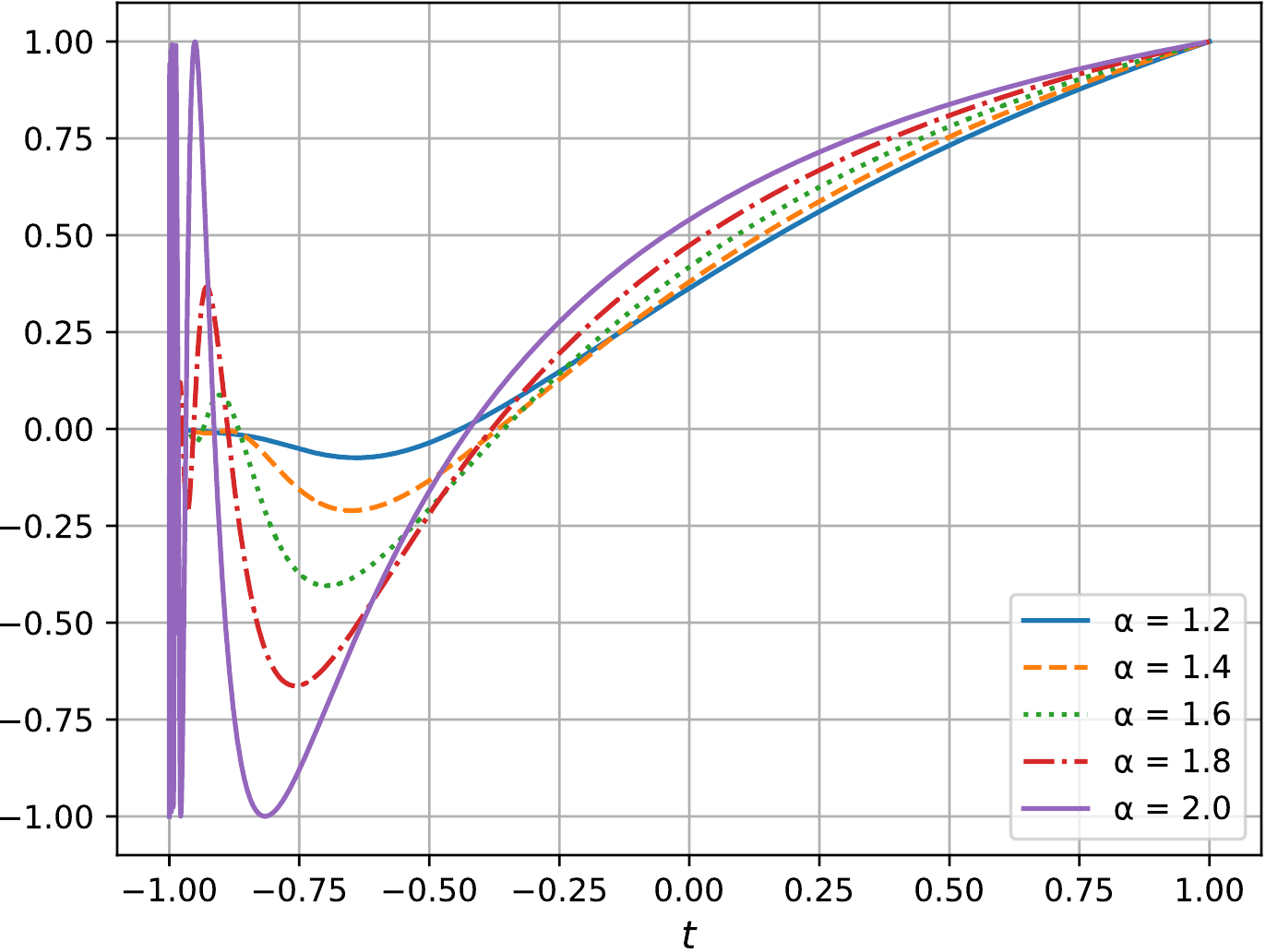} 
\end{center}
\caption{The function $f_\alpha(t)$ correponding to~$E_\alpha(-x)$
under the change of variable~\eqref{eq: change of variable}, for values
of~$\alpha$ in the range~$0<\alpha\le1$ (left) and $1<\alpha\le2$
(right).}\label{fig: f alpha}
\end{figure}

The substitution
\begin{equation}\label{eq: change of variable}
x=G(t)=\frac{1-t}{1+t}\quad\text{for $-1\le t\le1$}
\end{equation}
defines a one-one mapping $G:[-1,1]\to[0,\infty]$ with inverse given by
\[
t=G^{-1}(x)=\frac{1-x}{1+x}\quad\text{for $0\le x\le\infty$,}
\]
if we agree that $G(-1)=\infty$. Also, $G$ induces a 
bijection~$\mathcal{R}^m_m\bigl((-\infty,0])\bigr)\to 
\mathcal{R}^m_m\bigl((-1,1]\bigr)$ given 
by~$f\mapsto f_\sharp$, where $f_\sharp(t)=f\bigl(-G(t)\bigr)$.  Thus, it 
suffices to find the best approximation from~$\mathcal{R}^m_m([-1,1])$ to the 
function~$f_\alpha=(E_\alpha)_\sharp$, or in other words, the function
\begin{equation}\label{eq: f alpha}
f_\alpha(t)=E_\alpha\bigl(-x).
\end{equation}

\Cref{fig: f alpha} plots $f_\alpha$ for selected values of~$\alpha$ in 
the range~$0<\alpha\le2$. Recall from~\eqref{eq: integer alpha} that 
$E_0(-x)=1/(1+x)$ is already a rational function 
in~$\mathcal{R}^0_1\subseteq\mathcal{R}^1_1$.  If $0<\alpha\le1$, then 
$f_\alpha$ is well behaved and we can hope to approximate  it accurately by a 
function in~$\mathcal{R}^m_m([-1,1])$ with moderate values of~$m$. However, such 
an approximation becomes less and less feasible as~$\alpha>1$ increases because 
$f_\alpha(t)$ oscillates rapidly near~$t=-1$ due to the cosine factor in the 
term~\eqref{eq: cosine term}.

\begin{figure}
\begin{center}
\includegraphics[scale=0.4]{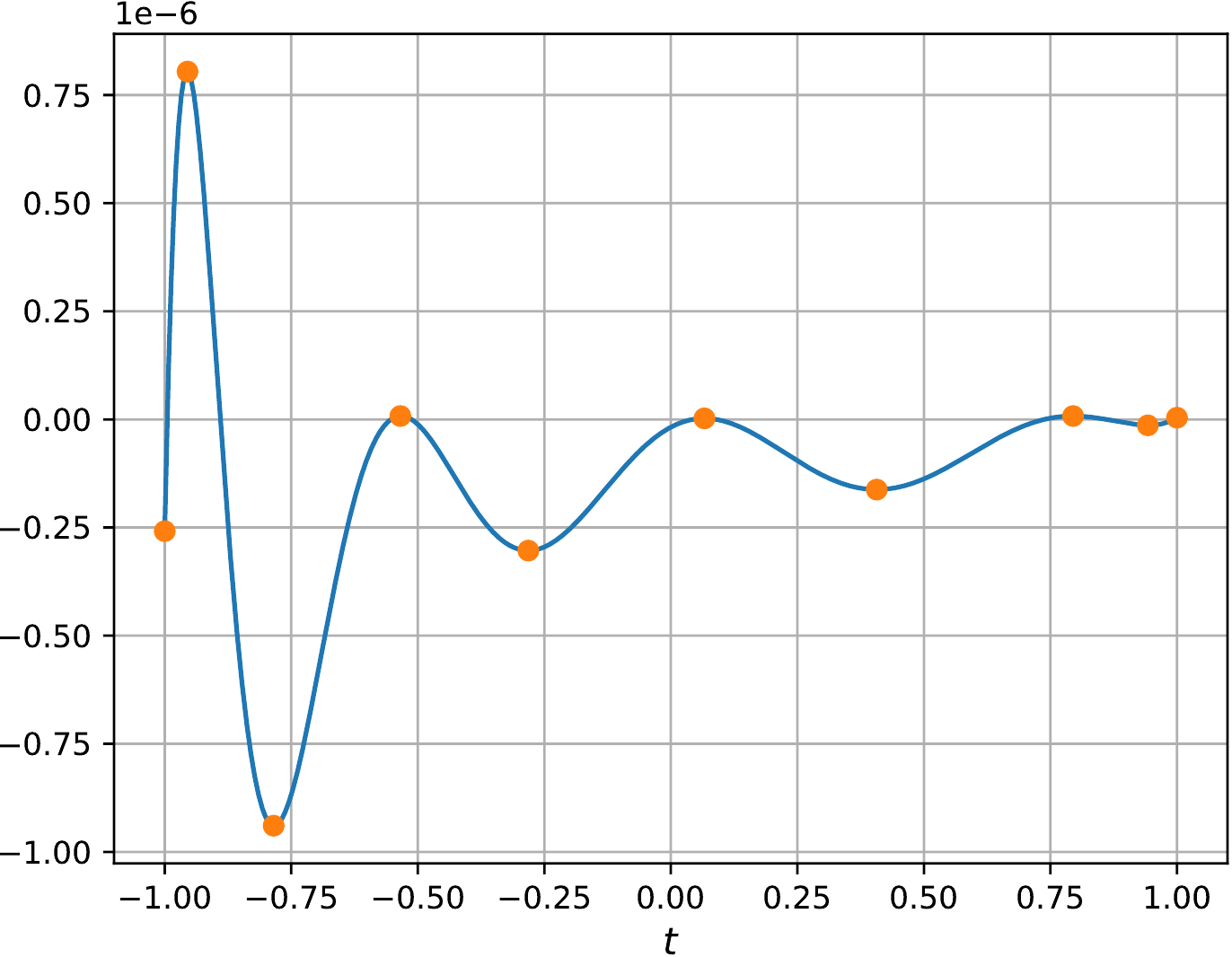} 
\hspace{0.4cm}
\includegraphics[scale=0.4]{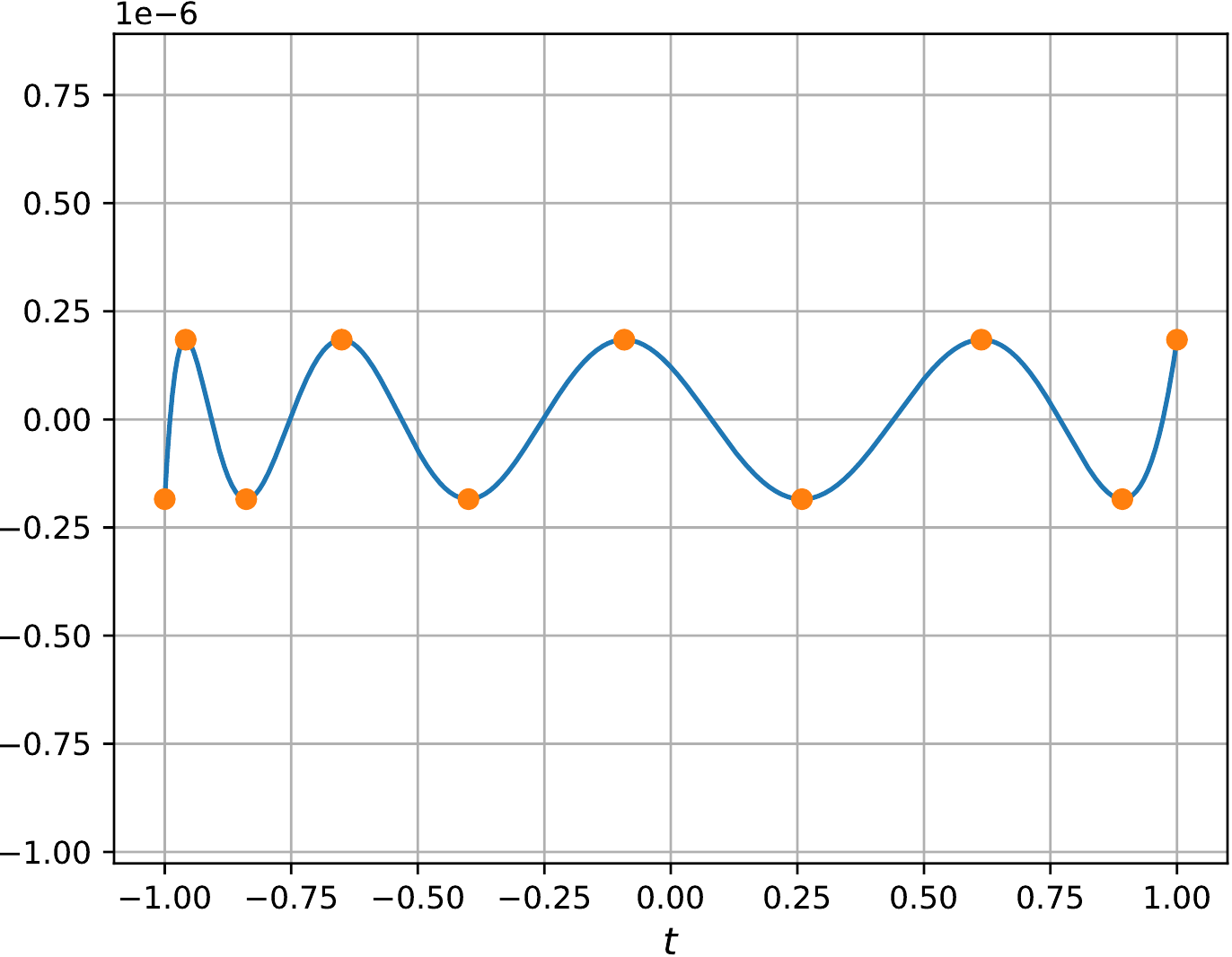} 
\end{center}
\caption{Errors in rational approximations from~$\mathcal{R}^m_m$ 
to~$f_{1/2}(t)=E_{1/2}(-x)$ for~$m=4$: (left) by the AAA algorithm and
(right) following subsequent iterations of the Remez algorithm.  The local
extrema are marked with dots.}
\label{fig: aaa remez} 
\end{figure}

\begin{table}
\caption{The uniform error $\|f_\alpha-R\|_{L_\infty(-1,1)}$ in the rational 
function~$R\in\mathcal{R}^m_m$ generated by the AAA algorithm.}
\label{tab: aaa errors}
\begin{center}
\renewcommand{\arraystretch}{1.2}
{\tt
\begin{tabular}{c|ccccc}
$m$&$\alpha=0{\cdot}2$&$\alpha=0{\cdot}4$&$\alpha=0{\cdot}6$&
$\alpha=0{\cdot}8$&$\alpha=1{\cdot}0$\\
\hline
 1 &   9.94e-01&   9.94e-01&   9.94e-01&   9.94e-01&   9.95e-01\\
 2 &   3.12e-03&   1.30e-02&   3.11e-02&   6.13e-02&   1.01e-01\\
 3 &   2.52e-05&   4.27e-04&   2.35e-03&   8.29e-03&   2.10e-02\\
 4 &   7.19e-08&   4.29e-06&   4.45e-05&   3.45e-04&   3.72e-03\\
 5 &   2.86e-10&   1.25e-07&   5.05e-06&   5.78e-05&   3.12e-04\\
 6 &   1.76e-12&   1.10e-09&   9.90e-08&   2.74e-06&   6.28e-05\\
 7 &   1.45e-14&   4.69e-11&   4.91e-09&   1.49e-07&   2.71e-06\\
 8 &   6.34e-13&   4.61e-13&   1.30e-10&   7.84e-09&   3.04e-07\\
 9 &   1.28e-12&   1.53e-14&   8.86e-12&   1.66e-09&   7.11e-08\\
10 &   4.30e-13&   1.90e-11&   1.83e-13&   5.49e-11&   5.17e-09
\end{tabular}
}
\end{center}
\end{table}

\begin{table}
\caption{The uniform error $\|f_\alpha-R^*\|_{L_\infty(-1,1)}$ in the best 
rational approximation~$R^*\in\mathcal{R}^m_m$.}\label{tab: best approx}
\begin{center}
\renewcommand{\arraystretch}{1.2}
{\tt
\begin{tabular}{c|ccccc}
$m$&$\alpha=0{\cdot}2$&$\alpha=0{\cdot}4$&$\alpha=0{\cdot}6$&
$\alpha=0{\cdot}8$&$\alpha=1{\cdot}0$\\
\hline
    1 &   5.00e-01&   5.00e-01&   5.00e-01&   5.00e-01&   5.00e-01\\
    2 &   2.08e-03&   8.55e-03&   2.02e-02&   3.87e-02&   6.68e-02\\
    3 &   7.72e-06&   1.30e-04&   7.13e-04&   2.55e-03&   7.36e-03\\
    4 &   2.83e-08&   1.94e-06&   2.48e-05&   1.66e-04&   7.99e-04\\
    5 &   1.04e-10&   2.89e-08&   8.61e-07&   1.07e-05&   8.65e-05\\
    6 &   3.80e-13&   4.31e-10&   2.98e-08&   6.93e-07&   9.35e-06\\
    7 &           &   6.41e-12&   1.03e-09&   4.47e-08&   1.01e-06\\
    8 &           &           &   3.56e-11&   2.89e-09&   1.09e-07\\
    9 &           &           &   1.23e-12&   1.35e-09&   1.17e-08\\
   10 &           &           &           &   1.20e-11&   1.26e-09\\
\end{tabular}
}
\end{center}
\end{table}

\begin{figure}
\begin{center}
\includegraphics[scale=0.6]{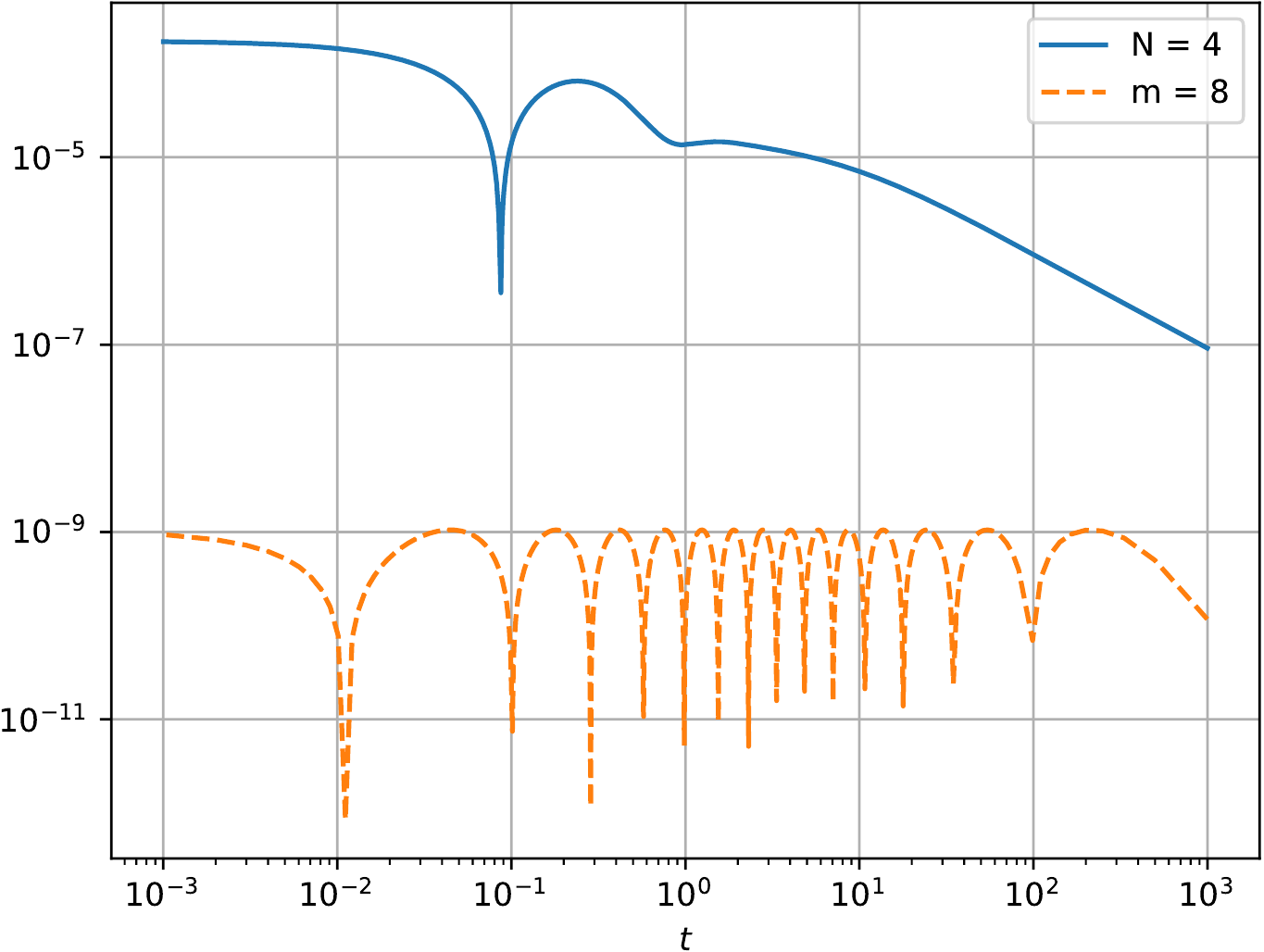} 
\end{center}
\caption{Solid line: the absolute error in the rational approximation 
$E_{\alpha,\beta}(-x)\approx Q_{\star,N}^{\mathrm{hyp}}(-x)$ from
$\mathcal{R}^{2N}_{2N+1}\bigl([0,\infty)\bigr)$ when $\alpha=0{\cdot}75$, 
$\beta=1$ and $N=4$.  Dashed line: the absolute error in the best approximation 
from~$\mathcal{R}^m_m\bigl([0,\infty)\bigr)$ for $m=2N=8$.}
\label{fig: compare rational} 
\end{figure}

For $0<\alpha\le1$, we applied the adaptive Antoulas--Anderson (AAA) 
algorithm~\cite{NakatsukasaEtAl2018} to generate an initial 
approximation~$f_\alpha\approx R\in\mathcal{R}^m_m([-1,1])$, and computed 
the local extrema $-1\le t_0<t_1<\cdots<t_{2m+1}\le1$ of the 
error~$f_\alpha-R$.  These points where then used to initialise a rational 
Remez algorithm~\cite{FilipEtAl2018} that generated the minimax 
approximation~$R^*\in\mathcal{R}^m_m([-1,1])$, which minimizes the error 
in the uniform norm and is characterized by the equioscillation property
\[
f_\alpha(t_l^*)-R^*(t_l^*)=(-1)^{l+1}\lambda\quad\text{for $0\le l\le2m+1$,}
\]
where $t_0^*$, $t_1^*$, \dots, $t_{2m+1}^*$ are the local extrema 
of~$f_\alpha-R^*$ and so $|\lambda|=\|f_\alpha-R^*\|_{L_\infty(-1,1)}$.  
\Cref{fig: aaa remez} shows the errors in these approximations 
when~$\alpha=1/2$ and~$m=4$.  In \cref{tab: aaa errors,tab: best approx}, we see 
how the uniform errors in $R$~and $R^*$ tend to zero as the polynomial 
degree~$m$ increases.  The convergence is most rapid for smaller values 
of~$\alpha$, as might be expected from the behaviour of~$f_\alpha$ seen in 
\cref{fig: f alpha}. The effects of roundoff become apparent once the errors 
reach about~$10^{-13}$, and the missing entries of \cref{tab: best approx} are 
cases when the Remez algorithm failed.  When computing $R$~and $R^*$, we 
evaluated $f_\alpha(t)=E_\alpha(-x)$ using the method of \cref{sec: real 
line} 
with~$N=14$ and hyperbolic contours.

Finally, recall from \cref{remark: rational} that the quadrature method 
generates a rational approximation~$E_{\alpha,\beta}(-x)\approx 
Q_{\star,N}(f;-x)$ from~$\mathcal{R}^{2N}_{2N+1}\bigl([0,\infty)\bigr)$
when $0<\alpha<1$.  \cref{fig: compare rational} compares the absolute error 
using this approximation with that of the best approximation 
from~$\mathcal{R}^m_m\bigl([0,\infty)\bigr)$ when $m=2N$ in the case 
$\alpha=0{\cdot}6$, $\beta=1$~and $N=4$.  The latter is smaller by $2$~to $5$ 
orders of magnitude.
\section{Conclusion}
The quadrature-based approach of \cref{sec: Quadrature} provides a practical 
method for numerical evaluation of the Mittag-Leffler 
function~$E_{\alpha,\beta}(z)$.  The maximum error can be reduced to around 
$10^{-N}$ using about $N$~terms, up to around $N=14$ when using standard 64-bit 
floating-point arithmetic; higher accuracy is achievable with larger~$N$ 
if extended precision is used. Somewhat better efficiency is possible for 
sufficiently small~$|z|$ using just the Taylor expansion~\eqref{eq: E defn}, 
and for sufficiently large~$|z|$ using the asymptotic 
expansion~\eqref{eq: asymptotics}. In the practically-important case when 
$0<\alpha<1$, rational approximation of~$E_{\alpha,\beta}(-x)$ 
for~$0\le x<\infty$ is effective.  
\bibliographystyle{spmpsci}
\bibliography{MLF_McLean}
\end{document}